\documentclass[11pt]{article}
\usepackage{GAFA2009}
\usepackage[paperwidth=193mm,paperheight=260mm,tmargin=25mm,bmargin=25mm,lmargin=25mm,rmargin=28mm]{geometry}
\usepackage{enumerate}
\usepackage{amscd}
\usepackage[mathscr]{eucal}
\usepackage{amsxtra}
\usepackage[all]{xy}
\begin{document}

\setfirstpage{0001}

\allowdisplaybreaks
     \abovedisplayskip=3.5pt plus 2pt minus 2pt
     \abovedisplayshortskip=-2pt plus 1pt minus 1pt
     \belowdisplayskip=3.5pt plus 2pt minus 2pt
     \belowdisplayshortskip=3.5pt plus 1pt minus 1pt

\newcommand{\TITLE}
{{
The {2}/{3}-convergence rate for the Poisson bracket
}}


\newcommand{\AUTHOR}{{L.\ Buhovsky}}

\newcommand{\AUtHOR}{{Lev Buhovsky}}

\Title
\AuthoR


\newcounter{theorem}[subsection]
\renewcommand{\thetheorem}{\thesubsection.\arabic{theorem}}

\newenvironment{theorem}[1][{}]{\refstepcounter{theorem}
      \smallbreak\noindent\begingroup%
      \textbf{Theorem
      \thetheorem{\rm{#1}}.}
      \endgroup\nobreak\slshape\ignorespaces}{\smallbreak}

\newenvironment{lem}[1][{}]{\refstepcounter{theorem}
      \smallbreak\noindent\begingroup%
      \textsc{Lemma
      \thetheorem{\rm{#1}}.}
      \endgroup\nobreak\slshape\ignorespaces}{\smallbreak}

\newenvironment{prop}[1][{}]{\refstepcounter{theorem}
      \smallbreak\noindent\begingroup%
      \textsc{Proposition
      \thetheorem{\rm{#1}}.}
      \endgroup\nobreak\slshape\ignorespaces}{\smallbreak}

\newenvironment{defn}[1][{}]{\refstepcounter{theorem}
      \smallbreak\noindent\begingroup%
      \textsc{Definition
      \thetheorem{\rm{#1}}.}
      \endgroup\nobreak\slshape\ignorespaces}{\smallbreak}

\newenvironment{que}[1][{}]{\refstepcounter{theorem}
      \smallbreak\noindent\begingroup%
      \textsc{Question
      \thetheorem{\rm{#1}}.}
      \endgroup\nobreak\rm\ignorespaces}{\smallbreak}

\newenvironment{exmp}[1][{}]{\refstepcounter{theorem}
      \smallbreak\noindent\begingroup%
      \textsc{Example
      \thetheorem{\rm{#1}}.}
      \endgroup\nobreak\rm\ignorespaces}{\smallbreak}

\newenvironment{remark}[1][{}]{\refstepcounter{theorem}
      \smallbreak\noindent\begingroup%
      \textsc{Remark
      \thetheorem{\rm{#1}}.}
      \endgroup\nobreak\rm\ignorespaces}{\smallbreak}


\def\Area{\operatorname{Area}}
\def\dist{\operatorname{dist}}
\Keywords{Symplectic manifold, Hamiltonian flow, Poisson
bracket,
rigidity, displacement energy, Hofer metric, uniform
norm, differential operator, Riemmanian metric}


\Thanks{This paper is part of the author's PhD thesis, being carried out
under the supervision of Prof.\ P.~Biran, at
Tel-Aviv University.
The author was partially supported by the Israel Science
Foundation (grant No. 1227/06 *)}

\vspace{-0.75cm}
\Abstract{
  In this paper we introduce a new method for approaching the
 $C^0$-rigidity results for the Poisson bracket.
Using this method, we
 provide a different proof for the lower semi-continuity under
 $ C^0 $ perturbations, for the uniform norm of the Poisson bracket. We
 find the precise rate for the modulus of the semi-continuity.
 This extends the previous results of Cardin--Viterbo,
Zapolsky,
 Entov and Polterovich. Using our method, we prove a
$C^0$-rigidity result
 in the spirit of the work of Humili\`{e}re. We also discuss a
 general question of the $C^0$-rigidity for multilinear
 differential operators.
}
\vspace{-0.5cm}

\section{Introduction and Main Results} \label{S:intro}

\subsection{Lower semi-continuity of the uniform norm of the Poisson
bracket.} \label{Lower}

The present note deals with the $C^0$-rigidity phenomenon
of
the Poisson bracket.
More precisely, for a symplectic manifold $(M,\omega)$,
we have a notion of a Poisson bracket $ \{\,\cdot\,{,}\,
\cdot\,\} : C^{\infty}(M) \times C^{\infty}(M) \rightarrow
C^{\infty}(M) $. For given $ f,g \in C^{\infty}(M) $ and a local
coordinate chart, this bilinear form involves partial derivatives
of the functions $ f,g $.
Therefore, we have no control of the
change of the values of $ \{ f,g \} $ when we perturb the
functions $ f,g  $ in the uniform norm. However, it turns out that
when we restrict ourselves to compactly supported functions on $ M
$, there exists a restriction on the uniform norm
$$
\big\| \{f,g \}
\big\|=\sup_{ x \in M }\big|\{ f,g \} (x)\big|\,,
$$
when we perturb $ f,g $
in the uniform norm. The first result in this direction was
obtained by F.\ Cardin and C.\ Viterbo \cite{CV}, who
showed that if $
\{f,g\}$ is not identically zero, then
$$
\liminf_{ \| F-f \|
,\,\| G-g \| \rightarrow 0 }\big\| \{ F,G \}\big\| > 0\,.
$$
This result
was improved by M.\ Entov, L.\ Polterovich, F.\ Zapolsky
(\cite{EPZ}, \cite{Z}, \cite{EP-1}).
It was shown in \cite{EP-1},
that in fact, for any symplectic manifold $ (M, \omega) $ and any
compactly supported $ f,g $, we have
$$
\liminf_{\|F-f\|,\,\|G-g\|\rightarrow
0}\big\|\{F,G\}\big\|
=\big\| \{ f,g \}\big\|\,.
$$
In both statements the functions $ F,G $ are
compactly supported.

\goodbreak
We introduce a new approach to the $C^0$-rigidity
phenomenon.
Our main result is summarized in Theorem
\ref{T:ExplicitBounds}.
Under the assumption that $\max\{f,g\}$ exists, we provide
an explicit lower estimate for the $ \sup \{ F,G \} $, when the
functions \hbox{$F,G:M\rightarrow\mathbb{R}$} are
$C^0$-close
to $ f,g $ respectively.

The statement of Theorem \ref{T:OpenManifold} coincides
with the abovementioned  result from \cite{EP-1}, while
stated under
slightly more general conditions.
In this case, our approach
enables us to provide a short proof of the statement.

In order to state the next theorem, we introduce the
following definition.

\begin{defn}
Let $(M,\omega)$ be a symplectic manifold. We denote by $
\mathcal{H}^{b}(M,\omega) $ the set of all smooth functions $ H :
M \rightarrow \mathbb{R} $, such that the Hamiltonian flow
generated by $ H $ is complete, that is, the solution exists for
any finite time.
\end{defn}

\begin{theorem} \label{T:OpenManifold}
\!Let $ (M,\omega) $ be a symplectic manifold.
Then, for any \hbox{$f,g\,{\in}\,C^{\infty}(M)$},
$$
\liminf_{F,G\in C^{\infty}(M),\,
 G \in \mathcal{H}^{b}(M,\omega),\,\|F-f\|,\|G-g\|
\rightarrow 0}\sup \{F, G \} = \sup \{ f,g \}\,.
$$
\end{theorem}

\vspace{3pt}
The method of the proof of Theorem \ref{T:OpenManifold}
is based on the positivity of the
 displacement energy of an open subset in $ M $(see  \cite{MS}).

\begin{defn}
Let $(M,\omega)$ be a symplectic manifold.
Given a pair of smooth functions $ f,g \in
C^{\infty}(M)$, we define
\begin{gather*}
  \Upsilon_{f,g}^{+}(\varepsilon):= \sup \{ f,g \}  - \inf_{ F,G \in C^{\infty}(M),
 G \in \mathcal{H}^{b}(M,\omega) ,  \| F-f \| \leqslant
 \varepsilon, \| G-g \|\leqslant\varepsilon}\sup\{F,G\}\,,
\\
\Upsilon_{f,g}(\varepsilon):=\big\|\{f,g\}\big\|-\inf_{F,
G \in C^{\infty}(M),
  G \in \mathcal{H}^{b}(M,\omega) ,  \| F-f \| \leqslant
\varepsilon,\|G-g\|\leqslant\varepsilon}\big\|\{F,G\}\big
\|\,.
\end{gather*}
\end{defn}

Then we have

\begin{theorem} \label{T:ExplicitBounds}
Let $ (M,\omega) $ be a symplectic manifold.
Assume that $f,g\in C^{\infty}(M)$ are such that $\{f,g\}
 $ attains its maximum at some $ x \in M $.
 Assume, in addition, that $ x $ is not a critical point
 for the functions $ f,g  $.
 Then
 $$
\limsup_{ \varepsilon \rightarrow 0 }
\frac{\Upsilon_{f,g}^{+}(\varepsilon)}{\varepsilon^{{2}/
{3}} } \leqslant 6
\big(-\{\{\{f,g\},f\},f\}(x)-  \{ \{ \{ f,g \} , g \} , g
 \}(x) \big)^{{1}/{3}} .
$$
Let us mention that, in the case of a closed manifold $
 (M,\omega) $, the condition that $ x $ is not a critical point for the functions
 $ f,g $ is satisfied automatically, if we
 assume that $ \{ f,g \} $ is not identically zero.
\end{theorem}

As will be seen from the proof of Theorem
\ref{T:ExplicitBounds}, the expression
 $$
-\big\{\{\{f,g\},f\},f\big\}(x)-\big\{\{\{f,g\} , g \} , g
\big \}(x)
$$
is non-negative, provided that the function $\{f,g\}$
attains its maximum
 at the \hbox{point $ x $}.

In the proof of Theorem \ref{T:ExplicitBounds} we use
lower
 estimates for the symplectic displacement energy.
 We use the notation $ e(W) $ for the symplectic displacement
 energy of the set $ W $.

For our purposes the following weak estimate will suffice.

\begin{prop}\label{P:SEnergy estimate}
Assume that we have a symplectic embedding
$$
i : U \subset (\mathbb{R}^{2n}, \omega_{std} )
 \hookrightarrow (M,\omega)\,.
$$
Consider a subset $V\subseteq U$ of the form
$V=Q_1\times Q_2\times\dots\times Q_n$, where\break
\hbox{$Q_1,Q_2,\ldots,Q_n\subset\mathbb{R}^{2}$}
are simply connected
planar
 domains.
Then we have
$$
e(i(V))\geqslant\tfrac{1}{2}\min\bigl(\Area(Q_1),\Area(Q_
2),\ldots
,\Area(Q_n)\big)\,.
$$
 \end{prop}

The Proposition \ref{P:SEnergy estimate} follows from the
 inequality (see \cite{MS})
$$
e(A)\geqslant\tfrac{1}{2}w_G(A)
$$
between the displacement energy $ e(A) $ of $ A $,
 and the Gromov width
$$
w_G(A)=\sup\big\{ \pi r^2 \mid B^{2n}(r) \text{ embeds
 symplectically in } A\big\}\,,
$$
where $ B^{2n}(r) \subset \mathbb{R}^{2n} $ is the standard
 Euclidean ball of radius $ r $.

  It is easy to see that replacing the functions $ f,F $ by $ -f, -F $ in
 Theorems \ref{T:OpenManifold} and \ref{T:ExplicitBounds},
we will get the analogous statements concerning the
$C^{0}$-rigidity of the
 infimum of the Poisson bracket. Both the rigidity of the supremum
 and of the infimum imply the corresponding rigidity result for
 the uniform norm $ \| \{ f,g \} \| $ of the Poisson bracket, since we have
$$
\big\|\{f,g\}\big\|=\max\Bigl(-\inf_{M}\{f,g\},\sup_{M}\{
f,
g\} \Big).
$$

   The coefficient $ 4 $ in the statement of the Theorem
\ref{T:ExplicitBounds}
 is not the exact value, and can be slightly improved using our method.
 On the other hand, weaker lower estimates of the form
 $$
e(i(V))\geqslant c\min\bigl(\Area(Q_1),\Area(Q_2),\ldots,
\Area(Q_n)\big)
 $$
for the displacement energy, will affect only this
coefficient, which will become larger.
 The precise optimal value is still to be found.

It turns out that the estimate on $ \Upsilon_{f,g}^{+}(\varepsilon)
$ in the Theorem \ref{T:ExplicitBounds} is sharp, up to
some constant factor.
 To obtain a lower bound for $ \Upsilon_{f,g}^{+}(\varepsilon)
 $, we first prove the following local result.

 \begin{theorem}  \label{T:LocalExample}
  Let $ (M,\omega) $ be a symplectic manifold.
Assume that we have $f,g\in C^{\infty}(M)$.
Denote by $ \Phi : M \rightarrow \mathbb{R} $
  the function
  $$
\Phi=-\big\{\{\{f,g\},f\},f\big\}-\big\{\{\{f,g\},g\} , g
 \big\}\,.
$$
Assume that $ \{ f , g \}
  $ attains its maximum at the point $ x \in M $, which is moreover a
  non-degenerate critical point of $ \{ f,g \} $.
Consider a neighborhood $ U $ of $ x $, and assume that
$$
\{f,g\}(y)<\{ f,g \} (x)\,,
$$
for every $ y \in \overline{U} \setminus \{ x \}
  $. Then we can find a neighborhood $ V  $ of $ x $, $
  \overline{V} \subset U $, such that for small $ \varepsilon > 0 $
  there exist smooth functions $ F,G : M \rightarrow \mathbb{R}
  $, satisfying
\begin{gather*}
\| F-f \| \leqslant
 \varepsilon\,,\quad\|G-g\|\leqslant\varepsilon\,,\\
\{ F,G \}(y)\leqslant\{f,g\}(x)-\tfrac{1}{3}\Phi(x)^{{1}/
{3} }
\varepsilon^{
{2}/{3}},\quad\forall y \in U\,,
\end{gather*}
 and such that $ F = f$, $G = g $
 on $ M \setminus V $.
 \end{theorem}

As a result of Theorems \ref{T:ExplicitBounds},
\ref{T:LocalExample}, we obtain
 the following global result on a closed manifold $ M $.

 \begin{theorem}  \label{T:GlobalEquivalence}
  Let $ (M,\omega) $ be a closed symplectic manifold.
Assume that we have $f,g\in C^{\infty}(M)$.
Denote by $ \Phi : M \rightarrow \mathbb{R}
  $ the function
  $$
\Phi=-\big\{ \{ \{ f,g \} , f \} ,f\big\}-\big\{\{\{f,g\}
, g \} , g
 \big\}\,.
$$
Assume that $x=x_{1},x_{2},\ldots,x_{N}$ are all the
points $ x \in M
$ for which $|\{f,g\}(x)|=\|\{f,g\}\|$, and assume that
all of them
  are non-degenerate critical points of the function $ \{ f,g \}
  $. Denote
$$
C=C(f,g)=\min\bigl(|\Phi(x_1)|,|\Phi(x_2)|,\ldots,|\Phi(x
_N)|\big)^{{1}/{3}}\,.
$$
  Then
$$
\frac{1}{3} C  \leqslant
  \liminf_{ \varepsilon \rightarrow 0 }\frac{\Upsilon_{f,
g}(\varepsilon)}{\varepsilon^{{2}/{3}}} \leqslant
  \limsup_{ \varepsilon \rightarrow 0 }\frac{\Upsilon_{f,
g}(\varepsilon)}{\varepsilon^{{2}/{3}}} \leqslant
6C\,.
$$
 \end{theorem}

It was shown in \cite{Z}, that in the case of dimension 2,
if $ \max_{M} \{ f,g \} $ is attained, then the statement of
 Theorem \ref{T:OpenManifold} in the dimension 2 case
becomes local
 in the sense of section \ref{S:Int and Loc} below, and
does not require the condition of $ G \in \mathcal{H}^{b}(M,\omega) $.
However, for dimensions bigger than 2, the situation
changes.
It turns out that the assumption $ G \in \mathcal{H}^{b}(M,\omega) $ in
Theorems \ref{T:OpenManifold}, \ref{T:ExplicitBounds} is
essential.
We show this in Example \ref{E:Gcondition} provided in
section \ref{S:Int and Loc}.
Moreover, Example \ref{E:NonLocal} in section
\ref{S:Int and Loc} shows
the non-locality of Theorem \ref{T:ExplicitBounds} for
any symplectic manifold $(M,\omega)$, with $\dim(M)>2$.
Examples \ref{E:Gcondition}, \ref{E:NonLocal} are closely
related, and we refer the reader to section
\ref{S:Int and Loc}
 for a detailed explanation of these phenomena.

After establishing the these results, the statement
of
Theorem \ref{T:ExplicitBounds} was re-proved by Entov and
Polterovich \cite{EP-2}, with the use of their own
approach.

\subsection{Conditions for the continuity of the Poisson
 bracket in the uniform norm.}
  Here we provide another application of the method, used to prove
 Theorems~\ref{T:OpenManifold}, \ref{T:ExplicitBounds}.
It is natural to ask the following:

\begin{que}  \label{Q:Unicontinuity}
Suppose we have a symplectic manifold $ (M,\omega) $,
functions \hbox{$f,g,h\in C^{\infty}(M)$}, and sequences
 $$
f_1,f_2,\ldots,g_1,g_2,\ldots\in C^{\infty}(M)\,,
$$
such
 that $ f_n \rightarrow f $, $  g_n  \rightarrow g $, $ \{ f_n ,g_n
 \} \rightarrow h $ uniformly on $ M $. Is
 it true that $ h = \{f,g \} $?
 \end{que}
 The answer in the general case is negative, as we see from the following example
 due to Polterovich.

 \begin{exmp} \label{E:PolterovichExample}
 On the plane $ \mathbb{R}^{2} $ consider the following sequence of
 functions:
$$
F_n(q,p)=\frac{\chi(p)}{\sqrt{n}}\cos(nq)\,,\quad G_n
(q,p)=\frac{\chi(p)}{\sqrt{n}}\sin(nq)\,,
$$
where $ \chi \in C^{
 \infty } ( \mathbb{R} ) $ given. Then $ \{ F_n ,G_n \} = \chi (p)
 \chi' (p) $, while $ F_n,G_n \rightarrow 0 $ uniformly.
 \end{exmp}

 We provide a sufficient condition under which we have an
 affirmative answer to this question.

Let us first introduce the notation needed for the
 formulation of the theorems in this section.

\begin{defn} \label{D:Norms1}
Suppose we have a smooth manifold $ X $ endowed with a
 Riemannian metric $ \rho $ and a smooth function $ h: X
 \rightarrow\mathbb{R} $. Take an integer $ k \geqslant 1 $.
For any $x\in X$, $v\in T_{x}X$, with the unit norm
$ \| v \|_{\rho} = 1 $, take a small $\rho$-geodesic
$ \gamma
 : [0, \varepsilon ) \rightarrow X $, such that $ \gamma (0)=x , \dot \gamma(0) = v
 $.
Then we denote
$$
\|h\|_{x,v,1}:=\bigg|\frac{d}{dt}|_{t=0} h(
 \gamma (t)) \bigg|.
$$
Next, for $ x \in X $, denote
$$
\| h \|_{x,1} := \max_{ v \in T_{x} X , \| v \|_{\rho} =
1}\|h\|_{x,v,1}\,.
$$
For a given subset $ Y \subset X $ with compact closure $
 \overline{Y} \subset X $, we denote
$$
\| h \|_{Y,1} := \sup_{ x \in Y } \| h \|_{x,1}\,.
$$
For a given subset $ Y \subset X $ with compact closure $
 \overline{Y} \subset X $, we denote
$$
\| h \|_{Y} :=
 \sup_{x \in Y } | h(x) |\,.
$$

We use the notation $\dist_{\rho}(x,y)$ for the
$ \rho $-distance between a pair of points
\hbox{$x,y\in X$}.
\end{defn}

 We first prove

 \begin{theorem} \label{T:UnicontinPrelim}
 Let $ (M,\omega) $ be a symplectic manifold, and an open subset $ U \subset M $ with compact closure
 $ \overline{U} \subset M $.
Assume that we are given a Riemmanian metric $ \rho $
 on $ U $, and smooth functions $ f,g \in C^{\infty}(M)
 $. Then there exists a constant $ C = C(U,\rho,f,g) > 0 $, such that
 for any $ F_1,G_1,F_2,G_2 \in C^{\infty}(M)$, satisfying
$$
\| F_1 -f \|_{U} ,\|F_2-f\|_{U},\| G_1 -g
 \|_{U} , \| G_2 -g \|_{U} < \varepsilon\,,
$$
we have
$$
\inf_{y,z\in U}\bigl|\{F_1,G_1\} (y) - \{ F_2 , G_2 \} (z)
\big|\leqslant C\varepsilon\max\bigl(1,\|G_1\|_{U,1},\|G_2
\|_{U,1} )\,.
$$
 \end{theorem}

As a corollary from Theorem \ref{T:UnicontinPrelim} we
obtain

 \begin{theorem} \label{T:UniContin}
 Let $ (M,\omega) $ be a symplectic manifold.
Assume that
 we have functions $ f,g,h \in C^{\infty}(M)$, and sequences
$$
f_1,f_2,\ldots,g_1,g_2,\ldots \in C^{\infty}(M)\,,
$$
such
 that $ f_n \rightarrow f $, $  g_n  \rightarrow g $, $ \{ f_n ,g_n
 \} \rightarrow h $ uniformly on $ M $.
Then if\break
\hbox{$\max(\|f_n-f\|_{U},\|g_n-g\|_{U})\|g_n\|_{U,1}
\rightarrow 0 $}
 for any open $ U \subset M $ with compact
 closure, then $ \{ f,g \} = h $.
The norms can be taken with respect to any Riemmanian
metric $ \rho $ \hbox{on $ M $}, and
 obviously the condition above does not depend on the metric.
 \end{theorem}

The proof of Theorem \ref{T:UniContin} uses Proposition
\ref{P:SEnergy estimate}.

 As it is easy to see, in Example \ref{E:PolterovichExample}
 we have
$$
\max\bigl(  \| F_n \| , \| G_n \|\big)\|G_n\|_{1}
\rightarrow \| \chi \|^{2}.
$$

The result of Theorem \ref{T:UniContin} is in the spirit
of the work of Humili\`{e}re \cite{H}.
Actually, he provides an affirmative answer to
Question \ref{Q:Unicontinuity}, if we assume that the
 sequences of pairs $ ( f_n, g_n ) $ of functions belong to some
 additional structure, namely a pseudo-representation of a normed Lie
 algebra.

  Using Theorem \ref{T:UnicontinPrelim}, one can extend the notion
 of Poisson bracket for some class of non-smooth functions.
 \begin{defn}
  Given a manifold $ X $, we say that the function $ f: X \rightarrow \mathbb{R} $ is of the H\"{o}lder class
 $ \alpha^+ $, if for some Riemmanian metric $ \rho $ on $ X $ and any $ x \in X $, we have
 $$
\lim_{\dist_{\rho}(x,y)\rightarrow 0}\frac{|f(x)-f(y)|}{(
\dist_{\rho}(x,y))^{\alpha}}=0\,.
$$
  Clearly the definition does not depend on the choice of the metric.
 \end{defn}

One can show that for given functions $ f,g : M
 \rightarrow \mathbb{R} $ of the H\"{o}lder class $ \frac{1}{2}^+
 $, one can define in a canonical way the analog of the Poisson bracket $ \{ f,g \}
 $, such that for any $ x \in M $, $ \{ f,g \}(x) $ is not a real number but a closed,
finite or infinite interval \hbox{in $\mathbb{R}$}.

 \subsection{Rigidity for general multi-linear
differential operators.}
 In this subsection we restrict ourselves to compactly
 supported functions.
 We ask the following general
 \begin{que}
 For a given smooth manifold $ X^n $, for which multi-linear
 differential operators on $ C^{\infty}(X) $, either of order 1 or bigger
 than 1, do we have some sort of $C^0$-rigidity?
 \end{que}

 We concentrate on the following two forms of $ C^{0} $ rigidity.

\begin{defn} \label{D:WeakAndStrongRigidity}
Assume that we have a multi-linear operator
$$
B: C^{\infty}(X)^{ \times m}
  \rightarrow C^{\infty}(X)\,.
$$
On the space $ C^{\infty}(X)^{ \times m} $
  consider the following metric:
given
$$
\mathcal{F}=(f_1,f_2,\ldots,f_m)\,,\quad\mathcal{G}=(g_1,
g_2,\ldots,g_m) \in
  C^{\infty}(X)^{ \times m},
$$
denote
$$
d_{C}(\mathcal{F},\mathcal{G}) := \max_{ 1
  \leqslant k \leqslant m} \| f_k - g_k \|\,.
$$
We say that $ B $ satisfies weak $ C^0 $-rigidity if,
for given $\mathcal{F}\in C^{\infty}(X)^{\times m}$, such
that $ \| B(
  \mathcal{F}) \| > 0 $ we have
$$
\liminf_{ d_{C}( \mathcal{ \widetilde{F} } ,
\mathcal{F})\rightarrow 0}\big\|B(\mathcal{\widetilde{F}})
\big\|>0\,.
$$
We say that $B$ satisfies strong $C^0$-rigidity if, for
given
  $ \mathcal{F} \in C^{\infty}(X)^{ \times m} $, we have
$$
\liminf_{d_{C}(\mathcal{\widetilde{F}},\mathcal{F})
\rightarrow 0 }
\bigl\|B(\mathcal{\widetilde{F}})\big\| = \| B(
  \mathcal{F}) \|\,.
$$
 \end{defn}

  On one hand, in the case of linear differential operators of the first order,
 the $ C^0 $-rigidity holds for
 any such operator, and moreover, it is local. We find an upper bound for the error,
 and it can be easily shown that it is precise, up to a constant factor.
 On the other hand, if we consider bilinear
 differential operators of the first order, then the necessary condition for $ C^{0} $
 rigidity is the anti-symmetricity of this form.
These statements are
 summarized in the following:

\begin{theorem}  \label{T:FirstOrder}
Consider a smooth manifold $ X^{n} $.
\begin{enumerate}[\indent\rm(a)]
\item
Suppose we are given a differential operator of the first order
$$
\lambda : C^{\infty}(X) \rightarrow C^{\infty}(X)\,,
$$
and a smooth function $ f: X \rightarrow \mathbb{R} $.
Assume that $ \lambda (f) $
attains its maximum at a point $ x $, such that $ x $ is a
non-degenerate critical point of $ \lambda (f) $. Take an
arbitrary open neighborhood $ U \subset X $ of $ x $. Then, for
any smooth function $F:X\rightarrow\mathbb{R}$ satisfying
$
\| F - f \|_{U} \leqslant \varepsilon $, we have
$$
\sup_{U}\lambda (F) \geqslant  \lambda (f)(x)  -\left(
\frac{9}{2}\right)^{{1}/{3}}\big(-\lambda^{3}(f)(x)\big)
^{{1}/{3}}\varepsilon^{{2}/{3} } - O(\varepsilon)\,.
$$
\item
Consider a bilinear differential operator of the first
order
$$
B(\,\cdot\,{,}\,\cdot\,):C^{\infty}(X)\times C^{\infty}(X
) \rightarrow C^{\infty}(X)\,,
$$
which is not antisymmetric. Then there exists a
function $ h \in C^{\infty}(X) $, and sequences $ f_n ,g_n \in
C^{\infty}(X) $ with $ \|f_n -h \|, \|g_n -h \| \rightarrow 0 $ ,
such that \hbox{$B(h,h)\neq 0$}, $B(f_n,g_n)=0$, for
every $n$.
\end{enumerate}
\end{theorem}

Let us focus on linear differential operators of the
first order.
First of all, the error is of the order $ \varepsilon^{
{2}/{3}}$, as we had in the case of the Poisson bracket.
This appears to be surprising because of the following
observation. Given a symplectic manifold $ (M,\omega) $,
and a function $ g \in C^{ \infty }(M) $, one can define the
linear operator $ \lambda (f) := \{ f,g \} $. On the other hand,
consider any differential operator of the first order on an
even-dimensional manifold $ X $. Then for any point $ x \in X $,
where the operator does not vanish, there exists a
neighborhood $
U $ of $ x $ and a symplectic structure $ \omega $ on $ U $, such
that our differential operator has the form $ \lambda (f) := \{
f,g \} $ on $ U $.

As we see, in Theorem \ref{T:ExplicitBounds}(a) we have
freedom in perturbing
both of the functions $ f,g $, while the application of
Theorem \ref{T:FirstOrder} allows us to perturb only one
of the
functions; nevertheless, this greater freedom does not decrease
the order of the error. Moreover, as an intermediate result in the
proof of Theorem \ref{T:ExplicitBounds}, we obtain
$$
\limsup_{
\varepsilon \rightarrow 0 }
  \frac{\Upsilon_{f,g}^{+}(\varepsilon)}{ \varepsilon^{2/3}
}\leqslant144^{{1}/{3}}
\Bigl(\max_{\theta}P(\theta)\Big)^{{1}/{3}}
,
$$
where $ P( \theta ) = - \{ \{  \{ f,g \} ,
  \cos( \theta )f + \sin( \theta
  )g \} ,  \cos( \theta )f + \sin( \theta )g \}(x) $.
 Replace the functions $ f,g $ by
$$
\cos( \theta )f + \sin( \theta )g  , - \sin( \theta )f +
\cos( \theta )g\,,
$$
for the value of $ \theta $, which gives us the maximum of $
P(\theta) $. Then the coefficient $ ( - \lambda^{3}(f) )^{
{1}/{3}}$ from Theorem \ref{T:FirstOrder} gives us the
exact
coefficient for the estimation of the error in
Theorem \ref{T:ExplicitBounds}, up to an absolute constant. Also
we see from the proof of Theorem \ref{T:LocalExample}, that in the
example which we provide there, we perturb only one of the
functions.

\begin{que} \label{Q:MultiOrder}
 Is it true, that in the case of general multi-linear
differential operators of the first order which satisfy the strong
version of $ C^0 $-rigidity, we also have this phenomenon? That
is, can the example which gives us the best error up to an
absolute constant be obtained by perturbing only one of the
functions?
\end{que}

As we see, the constant ${2}/{3} $ is not a special
symplectic constant.
We conjecture, that in fact the order $
\varepsilon^{{2}/{3}}$ for the error is correct for any
multi-linear differential operator of the first order,
which
satisfy the strong version of $ C^0 $-rigidity. It is evident
from the Theorem \ref{T:FirstOrder}, that it will be true,
provided the answer to Question
\ref{Q:MultiOrder} is affirmative.

Now we turn to the case of bi-linear differential
operators of the first
order.
It follows from Theorem \ref{T:FirstOrder} that in order
to have some $ C^{0}$-rigidity for a bilinear
differential operator
of the first order on $ C^{\infty}(X) $, it is necessary for this
operator to be anti-symmetric.  Actually, the statements of
Theorems \ref{T:OpenManifold}, \ref{T:ExplicitBounds}
show that
for a given manifold $ X $, their $ C^0 $-rigidity results hold
for all Poisson brackets derived from some given symplectic
structure $ \omega $ on $ X $, i.e.\ it holds for all
non-degenerate Poisson brackets on $ X $. However, taking an
arbitrary Poisson bracket on $ X $, not necessarily
non-degenerate, i.e.\ a bilinear operator
$$
\{\,\cdot\,{,}\,\cdot\,\} :
C^{\infty}(X) \times C^{\infty}(X) \rightarrow C^{\infty}
(X)\,,
$$
which is skew-symmetric, satisfies a Leibnitz rule and the Jacobi
identity, the manifold $ X $ is stratified into a disjoint union
of symplectic submanifolds, so we can reduce the situation to the
non-degenerate case.
Therefore, the statements of Theorems
\ref{T:OpenManifold}, \ref{T:ExplicitBounds} hold for any
Poisson structure on a smooth manifold $ X $. Observe that taking
a Poisson structure $ \{\,\cdot\,{,}\,\cdot\,\} $ on a
closed manifold $
X $, and a non-vanishing smooth function $ H(x) \in C^{\infty}(X)
$, we can define a new bilinear operator $ B(f,g) = H \cdot \{f,g
\} $. Then $ B $ will satisfy a weak form of $ C^{0} $ rigidity. A
priori, we cannot claim that $ B $ should satisfy the strong $ C^0
$-rigidity, because of the non-locality, presented in
Example \ref{E:NonLocal}.
However, if we assume that $ X $ admits
a fibration $ pr : X \rightarrow \mathcal{B}$ such that for any
fiber $ Y \subset X $, the values of $ \{ f,g \}|_{Y} $ depend
only on the restrictions $ f|_Y , g|_Y $, then, taking any
positive $ \mathcal{H} : \mathcal{B} \rightarrow \mathbb{R} $, the
form $ B(f,g)(x) = \mathcal{H}(pr(x)) \{ f,g \}(x) $ will satisfy
a strong form of rigidity, as can be easily seen.
For example, one
can take a 3-dimensional torus  $ \mathbb{T}^{3} = ( \mathbb{R} /
2 \pi \mathbb{Z} )^{3} $ with coordinates $ (x,y,z) \in
\mathbb{T}^{3} $, together with a fibration
$\mathbb{T}^{3}\rightarrow\mathbb{T}^{1}$,
$(x,y,z)\mapsto z$, and consider
$$
B(f,g)=\bigl(\sin(z)^2 + 1\big)(f_xg_y-f_yg_x)\,.
 $$
It is easy to
see that this particular $ B $ is not the Poisson bracket. As we
see, in this construction the form $ B $ is always degenerate.

\begin{que}
(a)
Is it true that, for closed manifolds the weak $ C^{0} $
rigidity holds only for multiples of a Poisson bracket by a
non-vanishing function?

(b) Is it true that for closed manifolds, in the case of
non-degenerate bilinear forms, the strong $ C^0 $-rigidity holds
only for Poisson brackets?

\end{que}

 Finally, the following example shows the existence of multi-linear operators of order 1, of
any number of functions, that satisfy the strong form of the $
C^{0} $-rigidity.

\begin{exmp}
 Given a natural $ m>1 $, take $ X =
 \mathbb{R}^{m} $, and define  $m$-linear $ \Phi :
 C^{\infty}(X)^{m} \rightarrow C^{\infty}(X) $ as follows: taking $
 f_1,f_2,\ldots,f_m \in C^{\infty}(X) $, define $ F: \mathbb{R}^{m}
 \rightarrow \mathbb{R}^{m} $ by $ F(x) := (
 f_{1}(x),f_{2}(x),\ldots,f_{m}(x)) $ and take $ B(f_1,f_2,\ldots,f_m) $
 to be the Jacobian $ J_F : \mathbb{R}^{m} \rightarrow \mathbb{R}
 $. The strong $ C^{0} $ rigidity for this $ B $ follows from
 simple volume considerations.
\end{exmp}

\subsection{Higher multiplicities of the critical points
of {\boldmath{$\{f,g\}$}}.}
Theorem \ref{T:ExplicitBounds}, applied to the case when
the
function $ \{ f,g \} $ has a degenerate maximum with multiplicity
bigger than 2 at the point $ x $, gives us only
$$
\Upsilon_{f,g}^{+}(\varepsilon)=o(\varepsilon^{{2}/{3}})
\,,
$$
without saying what is the order of $
\Upsilon_{f,g}^{+}(\varepsilon) $. It turns out that, after some
modification of the proof of Theorem
\ref{T:ExplicitBounds}, we
obtain

\begin{theorem} \label{T:Higher order}
 Let $ (M,\omega) $ be a symplectic manifold.

Assume that we have $f,g\in C^{\infty}(M)$, such that
$ \{ f , g \}
 $ attains its maximum at some $ x \in M $, and assume that the function $ \{ f,g \} $
 has multiplicity $ 2l $ at the \hbox{point $x$}.
 Assume in addition, that $ x $ is not a critical point
 for the functions $ f,g  $.
Define a differential operator
$$
\mathcal{D} (k) =
\big\{ \{ k , f \} ,f\big\} +\big\{ \{ k ,g\},g\big\}\,.
$$
 Then
$$ \limsup_{ \varepsilon \rightarrow 0 }
  \frac{\Upsilon_{f,g}^{+}(\varepsilon)}{
\varepsilon^{{2l}/({2l+1})}}\leqslant-9\left(\frac{1}{2l!}
\mathcal{D}^{l}(\{f,g\})(x)\right)^{{1}/({2l+1})}.
$$
\end{theorem}

The analogous statement holds also for the case of the
infimum.

\begin{remark}
Assume that $M$ is closed.
For every $ \varepsilon > 0 $, define a ``function"
\begin{gather*}
\mathcal{H}_{\varepsilon}
 : M \rightarrow \mathbb{R}\,,\\
\mathcal{H}_{\varepsilon} = \{ f,g \} + 9
\sum_{l=1}^{\infty} \varepsilon^{\frac{2l}{2l+1}} \left( \frac{1}{2l!}
\mathcal{D}^{l}(\{f,g\})\right)^{{1}/({2l+1})}.
\end{gather*}
Since this
series of functions does not have to converge, we consider $
\mathcal{H}_{\varepsilon} $ as a ``jet" in the functional space $
C^{\infty}(M) $, i.e.\ an asymptotic series, depending on
the
\hbox{parameter $ \varepsilon $}.
Then it is easy to see, that
Theorem \ref{T:Higher order} is equivalent to
$$
\inf_{ F,G \in C^{\infty}(M),\,
  G \in \mathcal{H}^{b}(M,\omega) , \, \| F-f \| \leqslant
\varepsilon,\,\|G-g\|\leqslant\varepsilon}\big\|\{F,G\}
\big\|
\geqslant
\| \mathcal{H}_{\varepsilon} \|\,,
$$
as ``jets". By this we mean
that for given $ L \geqslant 1 $, denoting the function $$
\mathcal{H}_{L,\varepsilon}=\{f,g\}+9\sum_{l=1}^{L}
\varepsilon^{\frac{2l}{2l+1}} \left( \frac{1}{2l!} \mathcal{D}^{l}( \{
f,g \} ) \right)^{{1}/({2l+1})},
$$
which is a truncation of the
asymptotic series $ \mathcal{H}_{\varepsilon} $, we have
$$
\inf_{ F,G \in C^{\infty}(M),\,
  G \in \mathcal{H}^{b}(M,\omega) , \, \| F-f \| \leqslant
\varepsilon,\,\|G-g\|\leqslant\varepsilon}\big\|\{F,G\}
\big\|
\geqslant
\| \mathcal{H}_{L,\varepsilon} \| - o\big(
\varepsilon^{{2L}/({2L+1})}\big)\,.
$$

In this observation, or reformulation of Theorem
\ref{T:Higher order},
we were able to collect all the cases of high multiplicities, and
moreover to get rid of considering all the critical points one by
one, and instead, to obtain a global inequality, which does not
apply to the critical points.
However, the asymptotic series $
\mathcal{H}_{\varepsilon} $ does not seem natural, because of the
possible non-smoothness of the functions, which enter in
its definition.
It would be interesting to find similar, but
correct, description of the result of Theorem
\ref{T:Higher order}.
Alternatively, it is possible that such a
description requires different framework and needs to be written
in other terms.
\end{remark}

\section{Proofs of Theorems} \label{S:proofs}

\PR{Proof of Theorem \ref{T:OpenManifold}}
 Let us first describe the main idea of the proof.

 We will use the notation $ X_{f},X_{g},X_{F},X_{G}  $ for the Hamiltonian vector fields
generated by the Hamiltonians $ f,g,F,G $ and by $
\Phi_{f}^{t},\Phi_{g}^{t}, \Phi_{F}^{t},\Phi_{G}^{t}  $ the
corresponding Hamiltonian flows.

 We have $ \{ f,g \} = df ( X_{g} ) $. Hence, roughly speaking, the
value of the Poisson bracket is the rate of change of values of
the function $ f $, computed through the Hamiltonian flow $
\Phi_{g}^{t} $ generated by $ g $. Assuming that, for some region
$ U \subset M $, we have $ \sup_{M} \{ F,G \} < \inf_{ U } \{ f,g \}
$, we will derive that for some small region $ W \subset U $ and
for some $ T > 0 $, the values of $ f(\Phi_{g}^{T}(W)) $ are
essentially bigger than those of $ F(\Phi_{G}^{T}(W)) $. If $ \|
F-f \| $ is small enough, the values of $ f(\Phi_{g}^{T}(W)) $
will be still much greater than those of $ f(\Phi_{G}^{T}(W)) $.
Hence, as a conclusion, we will get that the images $
\Phi_{g}^{T}(W) , \Phi_{G}^{T}(W) $ do not intersect, hence the
map $ \Phi_{g}^{-T} \circ \Phi_{G}^{T} $ displaces the set $ W $.
Using the positivity of the symplectic energy of $ W $, and the
upper estimate
$$
\| \Phi_{g}^{-T} \circ \Phi_{G}^{T} \|_{Hof}
\leqslant 2T \| g-G \| $$
on the Hofer norm, in the case when the
norm $ \| g-G \| $ is small enough, we will come to a
contradiction with our assumption that $ \sup_{M} \{ F,G \} < \inf_{
U } \{ f,g \} $.

 Let us turn now to the precise proof. Denote $ h= \{f,g \} $. Take any $ x \in M $ and
denote $ K = h (x) $.
Assume that, for some $ \delta > 0 $, we
have $ \{ F,G \} < K - \delta $ on $ M $, while $ \| f-F \|, \|
g-G \| < \varepsilon $. Here we will fix a specific $ \delta $,
while $ \varepsilon $ will be taken arbitrarily small.  For some
neighborhood $ U $ of $ x $, we will have that $ h(y) \geqslant K
- \frac{\delta}{2} $, for any $ y \in U $. Pick some $ V \subset U
$ and a positive $ T > 0 $, such that for any $ y \in V $, the
flow $ \Phi_{g}^{t}(y) $ exists for $ 0 \leqslant t \leqslant T $
and, moreover, $ \Phi_{g}^{t}(y) \in U $ for every $ 0 \leqslant t
\leqslant T $. Take an arbitrary point $ y \in V $ and define a
function $ K(t) = f( \Phi_{g}^{t}(y) ) $, $ t \in [0,T] $. Then we
have
$$
K'(t)=df\big( X_{g} ( \Phi_{g}^{t}(y) )\big) =\{f,g\}\big(
\Phi_{g}^{t}(y)\big)\geqslant K-\frac{\delta}{2}\,,
$$
for $ t \in [0,T] $.
Therefore, $ f( \Phi_{g}^{T}(y) ) - f(y) = K(T) - K(0)
\geqslant T\big(K-\frac{\delta}{2}\big)$.

 On the other hand, given any $ y \in M $, denote $ L(t) = F( \Phi_{G}^{t}(y) )
$, $ t \geqslant 0 $.
Then we have
$$
L'(t)=dF\big( X_{G} (\Phi_{G}^{t}(y))\big)=\{F,G\}\big(
\Phi_{G}^{t}(y)\big) \leqslant K -
\delta\,,
$$
for $ t \geqslant 0 $. Hence $ F( \Phi_{G}^{T}(y) ) -
F(y) = L(T) - L(0) \leqslant T(K - \delta) $. Since $ \| F-f \|
\leqslant \varepsilon $, we conclude that $ f( \Phi_{G}^{T}(y) ) -
f(y) \leqslant T(K - \delta) + 2 \varepsilon $.

Choose a small enough open subset $W\subset V$ such that
we have
$ | f(y) - f(z) | \leqslant{\delta T}/{3}$, when $y,z\in W
$. Then for any $ y,z \in W $ we have
\begin{align*}
f\big( \Phi_{G}^{T}(y)\big )&\geqslant T\left(K-\frac{
\delta}{2}\right)+f(y)\geqslant T\left(K -
\frac{\delta}{2}\right)-\frac{\delta T}{3}+f(z)\\
&\geqslant T\left(K-\frac{\delta}{2}\right)-\frac{\delta
T}{3} + f\big(
\Phi_{G}^{T}(z) \big) - T(K - \delta) - 2 \varepsilon\\
&=  f\big(
\Phi_{G}^{T}(z)\big )+\frac{\delta T}{6}-2\varepsilon\,.
\end{align*}
Assume
that $ \varepsilon <{\delta T}/{12} $.
Then we will get that
$f(\Phi_{G}^{T}(y))>f(\Phi_{G}^{T}(z))$ for any $ y,z \in
W $.
Therefore, $\Phi_{G}^{T}(W)\cap\Phi_{g}^{T}(W) = \emptyset
$, hence the map $ \Phi_{g}^{-T} \circ \Phi_{G}^{T} $ displaces
the set $ W $. Then, on one hand, the displacement energy $ e(W)
> 0 $, on the other hand we have an estimate for the Hofer norm:
$$
\|\Phi_{g}^{-T}\circ\Phi_{G}^{T}\|_{Hof}\leqslant 2T\|g-G
\| <
2T \varepsilon\,.
$$
Therefore, we conclude that $2T \varepsilon >
e(W) $. Observe that the choice of $ W , T $ depends only on $
f,g,x,\delta $.

As a conclusion, we get that, given $ f,g,\delta $, and
some point $ x \in M $,
there exists an open $ W \subset M $, and $ T>0 $, such that for
any $\varepsilon<\min({\delta T}/{12},{e(W)}/{2T})$
we have that for any $ F,G $ satisfying $ \| f-F \|, \| g-G \| <
\varepsilon $, we have $ \sup_{M} \{ F,G \} \geqslant \{ f,g \} (x)
- \delta $. Clearly this implies the statement of
Theorem \ref{T:OpenManifold}.
\qed

\PR{Proof of Theorem \ref{T:ExplicitBounds}}
The next definition describes the notation that will be
used
in the proof.

\begin{defn} \label{D:Norms2}
 Suppose we have a smooth manifold $ X $ endowed with a
 Riemannian metric $ \rho $ and a smooth function $ h: X
 \rightarrow\mathbb{R} $. Take an integer $ k \geqslant 1 $.
 For any $ x \in X $, $ v \in T_{x}X $
 with the unit norm $ \| v \|_{\rho} = 1 $, take a small $ \rho $-geodesic $ \gamma
 : [0, \varepsilon ) \rightarrow X $, such that $ \gamma (0)=x , \dot \gamma(0) = v
$.
Then we denote
$$
\|h\|_{x,v,k}:=\bigg|\frac{1}{k!}\frac{d^k}{dt^k}\Big|_{t
=0} h(
 \gamma (t))\bigg|.
$$
Next, for $ x \in X $ denote
$$
\|h\|_{x,k}:=\max_{v\in T_{x}X,\,\|v\|_{\rho}=1}\|h\|_{x,
v,k}\,.
$$
 For a given subset $ Y \subset X $ with compact closure $
\overline{Y}\subset X$, we denote
$$
\| h \|_{Y,k} := \sup_{ x \in Y } \| h \|_{x,k}\,.
$$
 For a given subset $ Y \subset X $ with compact closure $
\overline{Y}\subset X$, we denote
$$
\| h \|_{Y} :=
 \sup_{x \in Y } | h(x) |\,.
$$

Given a vector field $v$ on $X$, we denote by
$\|v\|_{x}=\| v(x) \| $ the norm of
 the vector $ v(x) \in T_{x} X $, with respect to $ \rho $. Then for a subset $ Y \subset X $
 with compact closure, we denote $ \| v \|_{Y} = \sup_{ x \in Y } \| v \|_{x} $.

We use the notation $\dist_{\rho}(x,y)$ for the
$\rho$-distance between a pair of points
\hbox{$x,y\in X$}.
\end{defn}

Note that for any $ Y \subset X $, $ \| \cdot \|_{Y,k} $
is not a norm, but rather a pseudo-norm on the space of smooth
functions.

We will use the notation $X_{f},X_{g},X_{F},X_{G}$ for
the Hamiltonian vector fields
generated by the Hamiltonians $ f,g,F,G $, and $ \Phi_{f}^{t},\Phi_{g}^{t},
\Phi_{F}^{t},\Phi_{G}^{t}  $ for the corresponding Hamiltonian flows.

The proof of Theorem \ref{T:ExplicitBounds} is a
generalization of the idea from the proof of
Theorem \ref{T:OpenManifold}. The proof can be divided into the
following parts. First, we consider functions $ f,g,F,G : M
\rightarrow\mathbb{R}$, such that
\begin{gather*}
\| f-F \| , \| g-G \| <
\varepsilon\,,\\
\max \{ F,G \} < \max \{ f,g \} - \delta\,.
\end{gather*}
We
take some neighborhood $ U $ of $ x $ in $ M $, and a Riemannian
metric $\rho $ on $ U $. We define some region $ W \subset U $,
depending on parameters $ \alpha,r $, and estimate the value range
of the function $ f $ on the images $ \Phi_{g}^{t}(W) ,
\Phi_{G}^{t}(W) $. We conclude that, under certain assumptions on
$ \varepsilon, \delta ,t $ and the parameters  $ \alpha, r $, the
images $ \Phi_{g}^{t}(W) , \Phi_{G}^{t}(W) $ do not intersect.
Therefore, under these assumptions, $W$ is displaced by
the map
$ \Phi_{g}^{-t} \circ \Phi_{G}^{t} $, hence we obtain
$$
2\varepsilon t >2\| g-G \| t \ge \| \Phi_{g}^{-t} \circ
\Phi_{G}^{t}\|_{Hof}\geqslant e(W)\,.
$$
On the other hand, we
find lower estimates for the displacement energy $ e(W) $ in
terms of $ \alpha, r $. Hence, under the assumptions on $
\varepsilon,\delta,t,\alpha,r$ above, and that
\begin{gather*}
 \| f-F \| , \| g-G \| <
\varepsilon\,,\\
\max \{ F,G \} < \max \{ f,g \} - \delta\,,
\end{gather*}
we
obtain an inequality concerning $ \varepsilon, \delta, t, \alpha,
r $.

 In the next step we consider $ f,g,F,G $, that satisfy $$ \| f-F
\| , \| g-G \| < \varepsilon ,$$ and we assume that we have such $
\delta,t,\alpha,r $, so that the abovementioned
assumption is
satisfied, but the inequality derived from the energy-capacity
argument is not.
Then we will have to conclude that
$$
\max \{ F,G \}\geqslant\max\{f,g\}-\delta\,.
$$
The next step in the
proof is to choose optimal $ t, \alpha, r $ to minimize $ \delta
$. The resulting formula involves estimations of $ C^2,C^1 $ norms
of $ \{f,g\} ,f,g $ on $ U $, with respect to the metric $ \rho $. Then we shrink the neighborhood $ U
$ to the point $ x $, arriving to the upper estimate for $ \delta
$, involving the norm of the Hessian of $ \{ f,g \} $, and norms
of $ X_{f},X_{g} $ at the point $ x $ with respect to the metric
$\rho$.

Finally, we choose the optimal metric $ \rho $ to obtain
 the statement of the Theorem~\ref{T:ExplicitBounds}.

  Let us turn to the proof. First of all, note that $ x $ is
  not a critical point for the functions $  f,g  $, and therefore
$$
df|_{x},dg|_{x},X_{f}(x),X_{g}(x) \neq 0\,.
$$
We start by choosing a Darboux neighborhood
$ i: U \hookrightarrow
(M,\omega)$ of $x$, where\break \hbox{$ 0 \in U \subset (
 \mathbb{R}^{2n},\omega_{std} ) $}, and $ i(0)=x $.
 Fix an arbitrary Riemannian metric $ \rho $ on $ i(U) $. Replacing $ U $ by some smaller open subset,
 we can guarantee that every point in $ i(U) $ can be joint to $ x $ by a $
 \rho $-geodesic, which lies in $ i(U) $.

 Then there exists an open neighborhood $ V \subset U $ of $ 0 $, and a positive $ T > 0 $, such that
 for any $ y \in i(V) $, the flow $ \Phi_{g}^{t}(y) $ exists when $ 0 \leqslant t \leqslant T
 $, and moreover, $ \Phi_{g}^{t}(y) \in i(U) $ for every $ 0 \leqslant
 t \leqslant T $.
Take some $0<r<\dist_{\rho}(x,M\backslash i(V))$ and some
real $
 \alpha > 0
 $, and consider the set
$$
W=W_{r,\alpha}=B_{x}(r)\cap\big\{y\in M\mid f(x)<f(y)<f(x)
+\alpha\big\} \subset M\,,
$$
where $ B_{x}(r) $ is a ball of radius $ r $ centered at $ x
 $, with respect to the metric $ \rho $.

 For $ y \in W $, denote $ K(t) = f( \Phi_{g}^{t}(y) ) $, $ t \in [0,T]
 $.
Then
$$
K'(t)=df\big(X_{g}(\Phi_{g}^{t}(y))\big)=\{f,g\}\big(
 \Phi_{g}^{t}(y)\big)\,.
$$
Denoting $ h = \{f,g \} $, we obtain that
$$
f\big(\Phi_{g}^{t}(y)\big)-f(y)=K(t)-K(0)=\int_{0}^{t}K'(
s) ds =
\int_{0}^{t}h\big(\Phi_{g}^{s}(y)\big)ds\,.
$$
 Let us estimate the value $ h(\Phi_{g}^{s}(y)) $ from below.
 First of all, we have
$$
\dist_{\rho}\big(x,\Phi_{g}^{s}(y)\big)\leqslant\dist_{
\rho}(x,y)
+\dist_{\rho}\big(y,\Phi_{g}^{s}(y)\big)\,.
$$
We have $\dist_{\rho}(y, \Phi_{g}^{s}(y) ) \leqslant
 s \| X_{g}\|_{U}$, $\dist_{\rho}(x,y) <
 r $, for $ y \in W $.
Hence $\dist_{\rho}(x, \Phi_{g}^{s}(y) )\break <
 r + s \| X_{g} \|_{U} $.

 \begin{lem} \label{L:HEstimate}
  For any $ z \in U $ we have
$$
h(z)\geqslant h(x)-\|h\|_{U,2}\dist_{\rho}
  (x,z)^{2} .$$
 \end{lem}

 \PR{Proof of Lemma \ref{L:HEstimate}}
  Take a $ \rho $- geodesic $ \gamma:[0,a] \rightarrow U $, such that
$$
\gamma(0)=x\,,\quad\gamma(a)=z\,,
$$
where $ a =\dist_{\rho}
  (x,z)$. Define $$ \phi : [0,a] \rightarrow
  \mathbb{R} $$ as $ \varphi (s):= h( \gamma (s)) $. Then, since
  the point $ x $ is a maximum point of $ h $, we have $ \varphi
  '(0) = 0 $. Therefore, $$ h(z) - h(x) = \varphi (a) - \varphi
  (0) = \int_{0}^{a} \varphi'(s) ds = \int_{0}^{a} ( a-s )
  \varphi''(s) ds\,.
$$
On the other hand, $ |\varphi''(s)| \leqslant 2 \| h
  \|_{U,2} $, so
$$
\big|h(z)-h(x)\big|\leqslant 2\|h\|_{U,2}\int_{0}^{a}(a-s
  ) ds = \| h \|_{U,2} a^{2} = \| h \|_{U,2}\dist_{\rho}
  (x,z)^{2} ,
$$
what implies the lemma.
  \qed

Hence for $ t \in [0,T] $ we have
  \begin{align*}
&f\big(\Phi_{g}^{t}(y)\big)-f(y)=
  \int_{0}^{t} h\big(\Phi_{g}^{s}(y)\big) ds\\
&\qquad>
\int_{0}^{t}h(x)-\|h\|_{U,2}\big(r+s\|X_{g}\|_{U}\big)^2
ds= h(x)t-\frac{1}{3}\frac{\|h \|_{U,2}}{\| X_{g} \|_{U}}
\big(r + t \| X_{g} \|_{U}\big)^3 ,
\end{align*}
so
  \begin{equation}  \label{e:1}
f\big(\Phi_{G}^{t}(y)\big)-f(y)>h(x)t-\frac{1}{3}\frac{\|
h \|_{U,2}}{\| X_{g} \|_{U}}
\big(r + t \| X_{g} \|_{U}\big)^3 .
  \end{equation}

Assume that we have smooth $F,G:M\rightarrow\mathbb{R}$
 and positive $ \varepsilon , \delta > 0 $, such that $$
 \parallel f-F \parallel , \parallel g-G \parallel < \varepsilon $$
 and
$$
\sup_{M}\{F,G\}<\max\{f,g\}-\delta=h(x)-\delta\,.
$$
Take some $z\in M$, and consider the function
$L(t)=F(\Phi_{G}^{t}(z)) $, $ t \geqslant 0 $.
We have
$$
L'(t)=dF\big(X_{G}(\Phi_{G}^{t}(z))\big)=\{F,G\}\big(
 \Phi_{G}^{t}(z)\big) < h(x) - \delta\,,
$$
hence we
 get an estimate
$$
L(t)-L(0)=F\big(\Phi_{G}^{t}(z)\big)-F(z)<\big(h(x)-
\delta\big)t\,,
$$
which holds for any $ z \in M $, $ t > 0 $.
 Since we have $ \parallel F-f \parallel < \varepsilon $, we
 obtain
\begin{equation}  \label{e:2}
f\big(\Phi_{G}^{t}(z)\big)-f(z)<\big(h(x)-\delta\big)t+2
\varepsilon\,.
 \end{equation}
In addition, for any $ y,z \in W $ we have
 \begin{equation}  \label{e:3}
 \big| f(y)-f(z) \big| < \alpha\,.
 \end{equation}

 From the inequalities (\ref{e:1}), (\ref{e:2}),
(\ref{e:3}) we derive,
that for any $y,z\in W$ we have
\begin{align*}
f\big(\Phi_{g}^{t}(y)\big)&>f(y)+h(x)t-\frac{1}{3}\frac{
\|h\|_{U,2}}{\|X_{g}\|_{U}}\big(r+t\|X_{g}\|_{U})^3\\
&> f(z) -
 \alpha + h(x)t - \frac{1}{3}\frac{\|h\|_{U,2}}{\|X_{g}\|
_{U}}\big(r + t \| X_{g} \|_{U}\big)^3\\
&>
f\big(\Phi_{G}^{t}(z)\big)-\big(h(x)-\delta\big)t-2
\varepsilon-\alpha+h(x)t-\frac{1}{3}\frac{\|h\|_{U,2}}{\|
X_{g} \|_{U}}\big(r + t \| X_{g}
\|_{U}\big)^3\\
&=f\big(\Phi_{G}^{t}(z)\big) + \delta t -\frac{1}{3}\frac
{\|h\|_{U,2}}{\|X_{g}\|_{U}}\big(r+t\|X_{g}\|_{U}\big)^3
 - 2\varepsilon - \alpha\,.
\end{align*}
 If we assume that
 \begin{equation}  \label{e:4}
\delta t\geqslant\frac{1}{3}\frac{\|h\|_{U,2}}{\|X_{g}\|_
{U}}\big(r + t \| X_{g} \|_{U}\big)^3
 + 2\varepsilon + \alpha
 \end{equation}
 holds, then for any $ y,z \in W $ we have
$$
f\big( \Phi_{g}^{t}(y)
 \big) > f\big(\Phi_{G}^{t}(z)\big)\,,
$$
therefore, the sets $ \Phi_{g}^{t}(W),
 \Phi_{G}^{t}(W) $ do not intersect.
Hence the map $ \Phi_{g}^{-t} \circ
\Phi_{G}^{t}$ \hbox{displaces $ W $}.
We have the following estimate
 for the Hofer norm:
$$
\|\Phi_{g}^{-t}\circ\Phi_{G}^{t}\|_{Hof}\leqslant 2t\|g-G
\|
 <2\varepsilon t\,.
$$
As a conclusion, we have the following:

 \begin{lem} \label{L:ConclusionExplBounds}
Assume now that we have smooth
$F,G:M\rightarrow\mathbb{R} $
  and positive $ \varepsilon , \delta > 0 $ such that $$
  \parallel f-F \parallel , \parallel g-G \parallel < \varepsilon $$
and
$$
\sup_{M}\{F,G\}<\max\{f,g\}-\delta=h(x)-\delta\,.
$$
In addition, assume that $(\ref{e:4})$ holds for some
$$
0 < t
\leqslant T\,,\quad 0<r<\dist_{\rho}(0,\partial V)\,,
\quad 0 < \alpha\,.
$$
Then for the set
$$
W=W_{r,\alpha}=B_{x}(r)\cap\big\{y\in M\mid f(x)<f(y)<f(x
) + \alpha\big\}\subset M\,,
$$
we have $2\varepsilon t > e(W) $.
 \end{lem}

Consider the case when we have smooth $ F,G : M \rightarrow \mathbb{R} $,
  positive $ \varepsilon , \delta > 0 $, and $ 0 < t
\leqslant T$, $0<r<\dist_{\rho}(0,\partial V)$,
$0<\alpha$, such
  that $ \parallel f-F \parallel , \parallel g-G \parallel < \varepsilon
$, the inequalities (\ref{e:4}) and
$2\varepsilon t\leqslant e(W) $ hold.
Then Lemma \ref{L:ConclusionExplBounds} will imply that
$$
\sup_{M} \{ F,G \} \geqslant
  \max \{ f,g \} - \delta\, .
$$

Assume that we have shown the existence of a positive
constant $
  C > 0$, such that if $ r, \alpha > 0 $ are small enough, and in
addition, ${\alpha}/{r} $ is small enough, then we have
$ e(W_{r, \alpha}) \geqslant C r
  \alpha $.
Then we will take $\alpha={2t \varepsilon
  }/{Cr} $, so that $2\varepsilon t \leqslant e(W) $.
  Then the inequality  (\ref{e:4}) is equivalent to
  \begin{equation}  \label{e:25}
    \delta  \geqslant \frac{\| h \|_{U,2}}{3} \frac{ (r + t \| X_{g}
   \|_{U})^3 }{ t \| X_{g} \|_{U} }
   + \frac{2\varepsilon}{t} + \frac{\alpha}{t} = \frac{\| h \|_{U,2}}{3} \frac{ (r + t \| X_{g}
   \|_{U})^3 }{ t \| X_{g} \|_{U} } + \frac{2\varepsilon}{t} +
\frac{2\varepsilon}{Cr}\,.
  \end{equation}

Our choice of $ t,r $ will be of the form
$t={P\varepsilon^{{1}/{3}}}/{ \| X_{g} \|_{U}} $,
  $ r = P \varepsilon^{{1}/{3}}$, for some $ P > 0 $.
Then we have
$$
\frac{\| h \|_{U,2}}{3} \frac{ (r + t \| X_{g}
  \|_{U})^3 }{ t \| X_{g} \|_{U} } + \frac{2\varepsilon}{t} +
\frac{2\varepsilon}{Cr} = \left(\frac{8}{3}\|h\|_{U,2}P^2
+2\left(\|X_{g}\|_{U}+\frac{1}{C}\right)\frac{1}{P}\right)
\varepsilon^{ \frac{2}{3} }.$$

Consider first the case, when $\|h\|_{U,2}>0$.
In this
case, the value of $ P $ that minimizes the expression
$$
\frac{8}{3} \| h \|_{U,2} P^2
+2\left(\|X_{g}\|_{U}+\frac{1}{C}\right)\frac{1}{P}\,,
$$
equals
$$
P = \bigg(
\frac{3}{8}\frac{\|X_{g}\|_{U}+\frac{1}{C}}{ \| h \|_{U,2}
}\bigg)^{{1}/{3}}.
$$
Then, for this $ P $,
$$
\frac{8}{3} \| h \|_{U,2} P^2
+2\left(\|X_{g}\|_{U}+\frac{1}{C}\right)\frac{1}{P}
=72^{{1}/{3}}\bigg(\|h\|_{U,2}
\left(\|X_{g}\|_{U}+\frac{1}{C}\right)^2\bigg)^{{1}/{3}}.
$$

In the case of $ \| h \|_{U,2} = 0 $, we fix arbitrary $ P > 0
  $.

  Note, that the choice of $ P $ we have made, does not depend on $
  \varepsilon $. We have
\begin{align*}
&t=\frac P{ \| X_{g} \|_{U}} \varepsilon^{ \frac{1}{3} }
\,,\\
&r=P\varepsilon^{{1}/{3} },\\
&\alpha  =
  \frac{2t \varepsilon }{Cr} = \frac{2}{C \| X_{g} \|_{U}}
\varepsilon\,,\\
&\frac{\alpha}{r} = \frac{2}{P C \| X_{g} \|_{U}
  } \varepsilon^{ \frac{2}{3} } .
\end{align*}

  Keeping the chosen value of $ P $ fixed, and taking $
  \varepsilon \rightarrow 0 $, we have
$$
t,\alpha,r,\frac{\alpha}{r}\rightarrow 0\,.
$$
In particular, $t\leqslant T$,
$r<\dist_{\rho}(0,\partial V)$, when $\varepsilon$ is
small enough. Moreover, for small
  enough $ \varepsilon $, the values of
$\alpha,r,{\alpha}/{r} $ are small, therefore we can
apply Lemma
\ref{L:SymplecticEnergy} to our situation.

\begin{lem} \label{L:SymplecticEnergy}
  For any $C<{1}/{\| X_{f} \|_{x}} $, we have
  $$
e( W_{r, \alpha } ) \geqslant C r \alpha\,,
$$
when $ \alpha , r,
{\alpha}/{r} \rightarrow 0 $.
  \end{lem}

  \PR{Proof of Lemma \ref{L:SymplecticEnergy}}
   We have $ W_{r, \alpha } \subset i(U) $, the Darboux
   neighborhood of $ x $. Take the pullback of $  W_{r, \alpha } $, the function $ f $ and the metric $ \rho $  to $
   U \subset ( \mathbb{R}^{2n} , \omega_{std} ) $, and denote the
   pullbacks by the same notation $ W_{r, \alpha }, f , \rho $.
Then in $U$ we have
$$
W_{r,\alpha}=B_{\rho,0}(r)\cap\big\{y\in\mathbb{R}^{2n}
\mid f(0) < f(y) < f(0) +
   \alpha\big\} \subset  \mathbb{R}^{2n} .
$$
Denote $b(\xi,\eta):=\rho|_{0}(\xi,\eta)$ the bilinear
form on
   $ \mathbb{R}^{2n} $, which is the restriction of $ \rho $ to
   the tangent space $ T_{0} ( \mathbb{R}^{2n} ) $. Denote $ l =
   df |_{0} $ - the differential of $ f $ at the point $ 0 $.
Then define
$$
\widetilde{W}_{r,\alpha} =
\big\{y\in\mathbb{R}^{2n}\mid b(y,y)<r^2\big\}\cap\big\{
y\in\mathbb{R}^{2n}\mid 0<l(y)<\alpha\big\}\subset
\mathbb{R}^{2n}.
$$
Then, for small $ r, \alpha $, we have $
   (1- o(1)) \widetilde{W}_{r,\alpha} \subseteq W_{r,\alpha} \subseteq
   (1+o(1)) \widetilde{W}_{r,\alpha} $.
Hence it is enough to
   establish $$ \frac{ e( i(\widetilde{W}_{r, \alpha } )) }{ r \alpha } \geqslant \frac{1}{\| X_{f} (0) \|_{\rho}}
-o(1)\,,
$$
when $r,\alpha,{\alpha}/{r}$ are small enough.
   Moreover, one can find a linear symplectic change of
   coordinates in $ \mathbb{R}^{2n} $, such that we will have $
l=df|_{0}=a\cdot dx_{1}$, for some $a\in\mathbb{R}$,
where $(x_{1},y_{1},\ldots
   , x_{n} ,y_{n} ) $ are coordinates in $ \mathbb{R}^{2n} $, so it is enough to consider this case only.
Denote $b_{11}=b({\partial}/{\partial
   y_{1}},{\partial}/{\partial y_{1}} ) $.
It is easy to see
   that for every $ 1 > \tau > 0 $, there exists some $ \kappa > 0 $, such that the set
$$
\big\{ y \in \mathbb{R}^{2n} \mid b(y,y) < r^2
 \big  \} $$
contains
$$
[-\kappa r,\kappa r]\times\left[-\tau\frac{r}{\sqrt{b_{11
}}},\tau\frac{r}{ \sqrt{b_{11}}}\right]\times[-\kappa r,
\kappa r]^{2n-2}\subset\mathbb{R}^{2n} \,,
$$
for any $ r > 0 $.
Hence the set $\widetilde{W}_{r,\alpha}$ contains
\begin{multline*}
[-\kappa r,\kappa r]\times\left[-\tau\frac{r}{\sqrt{b_{11
}}} , \tau \frac{r}{ \sqrt{b_{11}}}\right]\times[-\kappa
r,\kappa r]^{2n-2}\cap\left[0,\frac{\alpha}{a}\right]
   \times \mathbb{R}^{2n-2}\\
 =\left[0, \frac{\alpha}{a}\right] \times\left[-\tau\frac
{r}{\sqrt{b_{11}}},\tau\frac{r}{\sqrt{b_{11}}}
\right] \times [ -\kappa r, \kappa r]^{2n-2} ,
\end{multline*}
for small ${\alpha}/{r} $.
   We have that
$$
\Area\left(\left[0,\frac{\alpha}{a}\right]\times\left[-
\tau
\frac{r}{\sqrt{b_{11}}},\tau\frac{r}{\sqrt{b_{11}}}
\right]\right)=\frac{2\tau}{a\sqrt{b_{11}}}\alpha r\,,
$$
which is smaller than
$$
\Area\big([ -\kappa r, \kappa r] \times
    [ -\kappa r, \kappa r]\big) = 4 \kappa^2 r^2 ,
$$
when ${\alpha}/{r} $ is small enough.
Therefore, by Proposition \ref{P:SEnergy estimate}
   we have that the displacement energy
\begin{multline*}
e\left( i\left(\left [0, \frac{\alpha}{a} \right] \times
\left[-\tau\frac{r}{\sqrt{b_{11}}},\tau\frac{r}{\sqrt{b_{
11}}}\right]\times[-\kappa r,\kappa r]^{2n-2}\right)
\right)\\
\geqslant\frac{1}{2}\Area\left(\left[0,\frac{\alpha}{a}
\right]\times\left[-\tau\frac{r}{\sqrt{b_{11}}},\tau\frac
{r}{\sqrt{b_{11}}}
\right]\right)=\frac{\tau}{a\sqrt{b_{11}}}\alpha r\,.
\end{multline*}
Hence
$$
e\big(i(\widetilde{W}_{r,\alpha})\big)\geqslant e\left(
\left[0, \frac{\alpha}{a}\right] \times\left[-\tau\frac{r
}{\sqrt{b_{11}}},\tau\frac{r}{\sqrt{b_{11}}}
\right]\times[-\kappa r,\kappa r]^{2n-2}\right)  \geqslant
   \frac{ \tau}{ a \sqrt{b_{11}} } \alpha r\,.
$$
We have
$$ a^2 b_{11}
=a^2b\left(\frac{\partial}{\partial y_{1}},\frac{\partial
}{\partial y_{1}}\right)=b\left(a\frac{\partial}{\partial
y_{1}},a\frac{\partial}{\partial y_{1}}
\right)\,.
$$
Since $ df|_{0} = a \cdot dx_{1} $, then $
X_{f}(0)=a\frac{\partial}{\partial y_{1}}$, therefore
$$
a^2
b_{11}=b\left(a\frac{\partial}{\partial y_{1}},a\frac{
\partial}{\partial y_{1}}\right)=b\big(X_{f}(0),X_{f}(0
)\big)=\big\|X_{f}(0)\big\|_{\rho}^{2}\,,
$$
i.e.\ the square of the norm of the vector $X_{f}(0)$ with
   respect to the metric $ \rho $.
Therefore,
$$
e\big( i( \widetilde{W}_{r, \alpha } )\big )
   \geqslant  \frac{ \tau}{ a \sqrt{b_{11}} } \alpha r = \frac{ \tau}{ \| X_{f}(0) \|_{\rho} } \alpha
r\,,
$$
and this holds for any fixed $ 0 < \tau < 1 $, when we
   take $ \alpha, r $ to be small enough. This implies the lemma.
   \qed

Because of Lemma \ref{L:SymplecticEnergy}, we can take
arbitrary $ C <{1}/{\| X_{f} \|_{x}
  } $ . Then in the case of  $ \| h \|_{U,2} > 0 $, we can take
$$
\delta=72^{{1}/{3}}\bigg(\| h
  \|_{U,2} \left( \| X_{g} \|_{U} + \frac{1}{C}\right)^
{2} \bigg)^{{1}/{3} } \varepsilon^{{2}/{3} }
\,.
$$

  In the case of $ \| h \|_{U,2} = 0 $, for
  any fixed $ P > 0 $, we can take
$$
\delta=2\left(\|X_{g}\|_{U}+\frac{1}{C}\right)\frac{1}{P
} \varepsilon^{{2}/{3}}\,.
$$

Summarizing the above considerations, we see that if $ \| h \|_{U,2} > 0 $, then it follows
  that for any Darboux neighborhood $ i: U \hookrightarrow
  (M,\omega) $ of $ x $, and a Riemannian metric $ \rho $ on $ i(U)
$ we have
$$
\limsup_{\varepsilon\rightarrow 0}\frac{\Upsilon_{f,g}^{+}
(\varepsilon)}{ \varepsilon^{2/3}
}\leqslant 72^{{1}/{3}}
   \bigg( \| h \|_{U,2}  \left( \| X_{g} \|_{U} + \frac{1}
{C} \right)^2 \bigg)^{{1}/{3} } .
$$
Since this holds for any $C<{1}/{\|X_{f}\|_{x}}$, we
  obtain
$$
\limsup_{ \varepsilon \rightarrow 0 } \frac{\Upsilon_{f,g}^{+}(\varepsilon)}{ \varepsilon^{2/3}
}\leqslant 72^{{1}/{3}}
   \big( \| h \|_{U,2}  ( \| X_{g} \|_{U} + \| X_{f}
\|_{x})^2\big)^{{1}/{3} } .
$$

This inequality is correct also in the case of
$\|h\|_{U,2}=0
$, since then, fixing some specific $C<{1}/{\|X_{f}\|_{x}}
  $, we have
$$
\limsup_{ \varepsilon \rightarrow 0 } \frac{\Upsilon_{f,g}^{+}(\varepsilon)}{ \varepsilon^{2/3}
}\leqslant2\left( \| X_{g} \|_{U} + \frac{1}{C} \right)
\frac{1}{P}\,,
$$
  for any given $ P > 0 $, and hence $$ \limsup_{ \varepsilon \rightarrow 0 } \frac{\Upsilon_{f,g}^{+}(\varepsilon)}{ \varepsilon^{2/3}
  } = 0
$$
in this case.

Fixing the same metric $ \rho $ on $ U $, but shrinking $ U $
  to the point $ x $, we obtain
  \begin{equation}  \label{e:prelim}
   \limsup_{ \varepsilon\rightarrow 0}\frac{\Upsilon_{f,g}^{+}(\varepsilon)}{ \varepsilon^{2/3}
}\leqslant 72^{{1}/{3}}
\big(\|h\|_{x,2}(\|X_{g}\|_{x}+\|X_{f}\|_{x})^{2}\big)^{
{1}/{3}}.
  \end{equation}

The last step in the proof of the Theorem
\ref{T:ExplicitBounds} is to choose the optimal
\hbox{metric $\rho$} in
  the neighborhood of $ x $ in order to minimize the expression
on the right-hand side of the inequality
(\ref{e:prelim}).
From the inequality (\ref{e:prelim}) we see that it is
only essential to choose
  the metric on the tangent space $ T_{x}M $.

First consider the case when $X_{f}(x),X_{g}(x)\in T_{x}M$
  are linearly independent.
  In this case, the metric we choose will satisfy
  \begin{equation} \label{e:metric1}
\bigl\| \cos( \theta ) X_{f}  +
  \sin( \theta ) X_{g}\big \|_{\rho,x} = 1\,,
  \end{equation}
  for all  $ \theta $. It is easy to see that for any $ \varsigma > 0 $ we can find a
metric $\rho$ satisfying (\ref{e:metric1}), so that we
will
  have
  \begin{equation} \label{e:metric2}
\| h \|_{x,2} \leqslant \max_{\theta} \| h \|_{x, \cos( \theta )X_{f} +
   \sin( \theta )X_{g} ,2} + \varsigma\,.
  \end{equation}
  To do this, take any metric $ \rho $ which satisfies (\ref{e:metric1}), consider
  some linear complement of the linear subspace $ Sp(X_{f},X_{g}) \subset
  T_{x} M $, and then re-scale $ \rho $ by a sufficiently big
  factor in the direction of this complement.

   Assume now that we have a metric $ \rho $ that satisfies
  (\ref{e:metric1}), (\ref{e:metric2}). Suppose that for the
  vector $ v_{0} = \cos( \theta_{0} )X_{f} + \sin( \theta_{0}
  )X_{g}$ we have $$ \max_{\theta} \| h \|_{x, \cos( \theta )X_{f} +
\sin(\theta)X_{g},2}=\|h\|_{x,v_{0},2}\,.
$$
Then we have
$$
\|h\|_{x,2}\big(\| X_{g} \|_{x} + \| X_{f} \|_{x}\big)^{2}
  \leqslant 4 \| h \|_{x,v_{0},2} + 4 \varsigma\,.
$$

  We claim that
$$
\|h\|_{x,v_{0},2}=-\tfrac12\big\{\{h,\cos(\theta_{0})f+
\sin(
\theta_{0}
)g\},\cos(\theta_{0})f+\sin(\theta_{0})g\big\}(x)\,.
$$
In order to compute $\|h\|_{x,v_{0},2}$, we have to
choose a $ \rho $-geodesic $ \gamma
  : [0, \varepsilon ) \rightarrow M $, such that $ \gamma (0) = x , \dot \gamma(0) =
v_{0}$, and then
$$
\|h\|_{x,v_{0},2}=\left|\frac{1}{2}\frac{d^2}{dt^2}\Big
|_{t=0} h(
  \gamma (t))\right|.
$$
However, since $ h $ has at least order 2 at the
  point $ x $, we can only require from $ \gamma $ that $ \dot \gamma(0) =
  v_{0} $, without the assumption of being geodesic. In what
  follows, we can take $ \gamma(t) = \Phi_{k}^{t}(x) $, where $
  \Phi_{k}^{t} $ is the flow of the Hamiltonian $ k:=  \cos( \theta_{0} )f + \sin( \theta_{0} )g
  $. Then, denoting by $ X_{k} $ the Hamiltonian vector field of the Hamiltonian $ k $, we have
$$
\frac{d}{dt} h\big(\Phi_{k}^{t}(x)\big)=dh\big(X_{k}(
\Phi_{k}^{t}(x))\big)=\{h,k\}\big( \Phi_{k}^{t}(x) \big)
\,,
$$
hence
\begin{align*}
\frac{d^2}{dt^2}h\big(\Phi_{k}^{t}(x)\big)&=\frac{d}{dt}
\{h,k\}\big(\Phi_{k}^{t}(x)\big)=d\{h,k\}\big( X_{k}(
  \Phi_{k}^{t}(x))\big)\\
&=\big\{\{h,k\},k\big\}\big(\Phi_{k}^{t}(x)\big)
\,.
\end{align*}
Therefore, we have
\begin{align*}
\|h\|_{x,v_{0},2}&=\left|\left(\frac{1}{2}\frac{d^2}{dt
^2}\Big|_{t=0}h\big(\gamma(t)\big)\right)\right|=\left|
\frac{1}{2}\big\{\{h,k\},k\big\}(x)\right|\\
&=-\frac{1}{2}\big\{ \{  h,
  \cos( \theta_{0} )f + \sin( \theta_{0}
)g\},\cos(\theta_{0})f + \sin(\theta_{0})g\big\}(x)\,,
\end{align*}
  since $ x $ is the point of local maximum of $ h $.
   Hence we conclude that, denoting $ P( \theta ) = - \{ \{  h,
  \cos( \theta )f + \sin( \theta
  )g \} ,  \cos( \theta )f + \sin( \theta )g \}(x)  $, we
have
$$
\|h\|_{x,2}\big(\|X_{g}\|_{x}+\|X_{f}\|_{x}\big)^{2}
\leqslant 2\max_{\theta}P(\theta)+4\varsigma\,.
$$
So we have
\begin{align*}
\limsup_{ \varepsilon \rightarrow 0 }
\frac{\Upsilon_{f,g}^{+}(\varepsilon)}{ \varepsilon^{2/3}
}&\leqslant72^{{1}/{3}}
\Big(2\max_{\theta}P(\theta)+4\varsigma\Big)^{{1}/{3}}\\
&=144^{{1}/{3}}
\Big(\max_{\theta}P(\theta)+2\varsigma\Big)^{{1}/{3}}
\,.
\end{align*}
Since this holds for any $
  \varsigma > 0 $, we obtain
$$
\limsup_{ \varepsilon \rightarrow 0 }
  \frac{\Upsilon_{f,g}^{+}(\varepsilon)}{ \varepsilon^{2/3}
}\leqslant144^{{1}/{3}}
\Big( \max_{\theta} P( \theta )\Big)^{{1}/{3}}.
$$
It is easy to see that
$P(\theta)+P\big(\theta+\frac{\pi}{2}\big) =
  - \{ \{ h , f \} ,f \}(x) -  \{ \{ h , g \} , g \}(x) $ for
\hbox{every $\theta$}, and since $x$ is a local maximum
point of $
  h $, we have $ P(\theta) \geqslant 0 $ for every $ \theta $.
  This implies $ \max_{\theta} P(\theta) \leqslant  - \{ \{ h , f \} ,f \}(x) -
  \{ \{ h , g \} , g \}(x) $.
Therefore,
$$
\limsup_{ \varepsilon \rightarrow 0 }
  \frac{\Upsilon_{f,g}^{+}(\varepsilon)}{ \varepsilon^{2/3}
}\leqslant144^{{1}/{3} }
\big(  - \{ \{ h , f \} ,f \}(x) -
  \{ \{ h , g \} , g \}(x)\big)^{{1}/{3} } .
$$

It remains to check the case when
$X_{f}(x),X_{g}(x)\in T_{x} M $
  are linearly dependent. Suppose for instance that $ X_g = q X_f
  $, when $ | q | \leqslant 1 $ (the other case is similar).
  Take any metric $ \rho $, such that  $ \| X_{f}  \|_{\rho,x} =
  1 $, then take some $ \varsigma > 0 $, and re-scale $ \rho $ along some linear complement of $
  Span( X_{f} ) $, so that we will have
  \begin{equation} \label{e:metric3}
\| h \|_{x,2} \leqslant \| h \|_{x, X_{f} ,2}+\varsigma\,.
  \end{equation}
   We have
$$
\|h\|_{x,X_{f},2}=-\tfrac{1}{2}\big\{\{ h,f\},f\big\}(x)
\,,
$$
therefore
$$
\|h\|_{x,2}\big(\| X_{g} \|_{x} + \| X_{f} \|_{x}\big)^{2}
  \leqslant -2\big\{ \{  h,f\},f\big\}(x)+4\varsigma\,.
$$
Hence
\begin{align*}
\limsup_{ \varepsilon \rightarrow 0 }
  \frac{\Upsilon_{f,g}^{+}(\varepsilon)}{ \varepsilon^{2/3}
}&\leqslant72^{{1}/{3}}
\bigl(-2\{\{h,f\},f\}(x)+4\varsigma\big)^{{1}/{3}}\\
&=144^{{1}/{3}}
\bigl(-\{\{h,f\},f \} (x)  + 2 \varsigma\big)^{{1}/{3} }
\,.
\end{align*}
Since this holds for any $\varsigma>0$, we obtain
\begin{align*}
\limsup_{ \varepsilon \rightarrow 0 }
  \frac{\Upsilon_{f,g}^{+}(\varepsilon)}{ \varepsilon^{2/3}
}&\leqslant144^{{1}/{3}}
\bigl(  - \{ \{  h, f \} , f \} (x)\big)^{{1}/{3} }\\
&\leqslant144^{{1}/{3}}
\bigl(-\{\{h,f\},f\}(x)-\{\{h,g\},g\}(x)\big)^{{1}/{3} }
\, .
\end{align*}

Since $144^{{1}/{3}
  } <6$, we obtain the desired result.
\qed

\PR{Proof of Theorem \ref{T:LocalExample}}
  Denote by $ X_f , X_g $ the Hamiltonian vector fields generated by Hamiltonians $
 f,g: M \rightarrow \mathbb{R} $. Denote $ h= \{ f,g \} $. Since $ x $ is the
 local maximum point of
 $ h $, we have
$$
\big\{\{h,f\},f\big\}(x),\big\{\{h,g\},g\big\}(x)
\leqslant
  0\,.
$$
 If $ \{ \{ h,f \} , f \}(x) = \{ \{ h,g \} , g \}(x) = 0 $, there is nothing to prove.
 Consider the complementary case.
 Without loss of generality, we can assume that $ \{ \{ h,g \} , g \}(x)
 < 0 $, $ \{ \{ h,g \} , g \}(x)\leqslant\{\{h,f\},f\}(x)
$ (in the opposite case, we can apply the
 Theorem \ref{T:LocalExample} to the functions $ -g, f $).
  Because of $  \{ \{ h, g \} , g \}(x) < 0 $, we have
 $ X_{g}(x) \neq 0 $.
Hence, for some small neighborhood $ W
 \subset U $ of $ x $, there exists a coordinate $ x_{1} : W
\rightarrow\mathbb{R}$, such that $x_{1}(x)=0$,
$X_{g}={\partial}/{\partial x_1
 } $ on $ V $. Denote $ H = h_{x_{1}} $. Then $ H_{x_{1}} = \{ \{ h, g \} , g
 \} \neq 0 $, therefore one can extend $ x_{1} $ to a coordinate
 system $ (x_1, y_1,x_2,y_2,\ldots,x_n,y_n ) $ on $ W $, such that
$$
H_{y_1}(x)=H_{x_2}(x)=H_{y_2}(x)=\dots=H_{x_n}(x)=H_{ y_n
 }(x)= 0\,.
$$
Note that this is not necessarily a Darboux coordinate
system.
Denote
$A=\break-\{\{h,g\},g\}(x)= -h_{ x_1 x_1 } (x) > 0 $.
Take some $ b>0 $ , such that the cube
$$
K=\big\{
(x_1,y_1,x_2,y_2,\ldots,x_n,y_n)\mid-b\leqslant x_1,y_1,x
_2,y_2,\ldots,x_n,y_n\leqslant b\big\}
$$
is
 inside $ W $. Denote also
$$
K'=\big\{
 y=(x_1, y_1,x_2,y_2,\ldots,x_n,y_n
 ) \in K \mid-b/3\leqslant x_1\leqslant b/3\big\}\,.
$$
 For small $ \varepsilon > 0 $, take a smooth $ \varphi : \mathbb{R} \rightarrow
 \mathbb{R} $, such that
$\varphi(t)=\frac{1}{2}A^{{1}/{3} }
\varepsilon^{{2}/{3}}t$ for $t\in[-A^{-{1}/{3}
}\varepsilon^{{1}/{3}},A^{-{1}/{3} }
 \varepsilon^{{1}/{3}}]$,
that $ \varphi'(t) \geqslant 0 $
 when $ t \in [-b/3,b/3] $, that
$$
\varphi'(t) \geqslant \frac{
 \max_{y \in K \setminus K'} h(y) - h(x) }{2} $$ for $ t \in [-2b/3, 2b/3] $,  $ \varphi(t) = 0 $ for $ t \in [-b,-2b/3] \cup [ 2b/3,
 b] $ , and $ |\varphi(t)| \leqslant \varepsilon $ for any $ t \in \mathbb{R} $.
Then take some bump function
$\psi:\mathbb{R}^{2n-1}\rightarrow\mathbb{R}
 $, such that $ \psi = 1 $ on $ \frac{1}{3} K $ and $ \psi = 0 $ outside $ \frac{2}{3} K
 $, and $ 0 \leqslant \psi \leqslant 1 $ on $ \mathbb{R}^{2n} $.
Then define $ F,G: M \rightarrow \mathbb{R} $ by $ F = f $
 on $ M \setminus W $, and $ F = f - \varphi(x_1)
 \psi(y_1,x_2,y_2,\ldots,x_n,y_n) $ on $ W $, and then take $ G=g $
\hbox{on $ M $}.
Note that $ F=f $ on $ W \setminus K $.

First of all, for any
$y=(x_1,y_1,x_2,y_2,\ldots,x_n,y_n)\in W$, we have
$$
\big| f(y) - F(y)\big | =
\big | \varphi(x_1) \psi(y_1,x_2,y_2,\ldots,x_n,y_n)\big | \leqslant\big |
 \varphi(x_1)\big|\leqslant\varepsilon\,.
$$
For $ y \notin W $ we
 have $ f(y) -F(y) = 0 $.
Therefore, $ \| f-F \| \leqslant
 \varepsilon $. As $ G=g $, we have $ \| g-G \| = 0 \leqslant
 \varepsilon $.
  On the other hand, for any function $ k: W \rightarrow
\mathbb{R}$, we have
$$
\{ k,g \} = dk( X_{g} ) = dk \left( \frac{\partial}{\partial
 x_1} \right) = k_{x_1}\,.
$$
Therefore, for $ y = (x_1,y_1,x_2,y_2,\ldots,x_n,y_n) \in W $  $$ \{ F,G \} = \{ f
 - \varphi \psi , g \} = \{ f,g \} - \{ \psi\varphi ,g \} = h - \varphi'(x_1)
 \psi(y_1,x_2,y_2,\ldots,x_n,y_n)\,.
$$
We wish to show that
$\{F,G\}\leqslant\{f,g\}-\frac{1}{2}A^{{1}/{3}}
\varepsilon^{{2}/{3}}
$ on $ W $. This is equivalent to $$ \varphi'(x_1)
\psi(y_1,x_2,y_2,\ldots,x_n,y_n)\geqslant h(y)-h(x)+
\tfrac{1}{2}A^{{1}/{3} } \varepsilon^{
{2}/{3} } .
$$
Because of the condition
\begin{gather*}
h_{ x_1 x_1 }(x)=-A\,, \\
 h_{x_1 y_1}(x) = h_{ x_1 x_2 }(x) = h_{ x_1 y_2 }(x)
=\dots= h_{ x_1 x_n }(x) = h_{ x_1 y_n }(x)= 0\,,
\end{gather*}
and since $ x $ is a non-degenerate critical point of $ h $,
we have that the domain
$$
\left\{y\in W\bigm| h(x)-h(y)\leqslant\tfrac{1}{2} A^{
{1}/{3} } \varepsilon^{
{2}/{3} }\right\}
$$
lies inside the set
$$
K''=\big\{ y=(x_1,y_1,\ldots,x_n,y_n) \in W
\,,\ |x_1|\leqslant A^{-{1}/{3}}\varepsilon^{{1}/{3}
 }\big\}\cap\tfrac{1}{3}K\,,
$$
when $ \varepsilon $ is small.
 For $ y \in K'' $,
$$
\varphi'(x_1)\psi(y_1,x_2,y_2,\ldots,x_n,y_n)=\tfrac{1}{2
}A^{{1}/{3}} \varepsilon^{
{2}/{3} } \geqslant
 h(y)-h(x)+\tfrac{1}{2}A^{{1}/{3} } \varepsilon^{
{2}/{3}}.
$$
For $y\in K'\setminus K''$,
$$
\varphi'(x_1)\psi(y_1,x_2,y_2,\ldots,x_n,y_n)\geqslant 0
\geqslant h(y)-h(x)+\tfrac{1}{2} A^{{1}/{3}}\varepsilon^{
{2}/{3}}.
$$
For $y\in K\setminus K'$,
\begin{align*}
\varphi'(x_1)
 \psi(y_1,x_2,y_2,\ldots,x_n,y_n)&\geqslant \frac{
 \max_{z \in K\setminus K'}h(z)-h(x)}{2}\\
& \geqslant h(y) - h(x) +
 \tfrac{1}{2}A^{{1}/{3} } \varepsilon^{
{2}/{3}},
\end{align*}
when $ \varepsilon $ is small. Since $ F = f , G = g $ on $ U \setminus K $,
 and
$$
\sup_{y\in U\setminus K}h(y)<h(x)\,,
$$
we have
$$
\{ F,G \} (y) = \{ f,g \} (y) = h(y) \geqslant
 h(x) -  \tfrac{1}{2}A^{{1}/{3}}\varepsilon^{{2}/{3}}
$$
for $ y \in U \setminus K $,
 when $ \varepsilon $ is small.
Hence we have shown that for $ V:=int(K) \subset U $, for $ \varepsilon $ small enough, there exist
 smooth $ F,G : M \rightarrow \mathbb{R} $, such that $ F=f,g=G $ on $ M \setminus V $, and
\begin{gather*}
\| F-f \| \leqslant
\varepsilon\,,\quad\|G-g\|\leqslant\varepsilon\,,\\
\{F,G\}(y)\leqslant\{f,g\}(x)-\tfrac{1}{2}\big(-\{\{h,g\},
g\}(x)\big)^{{1}/{3} }
 \varepsilon^{{2}/{3}} ,\quad\forall y\in U\,.
\end{gather*}
We have
$$
\tfrac{1}{2}\big(-\{\{h,g\},g\}(x)\big)^{{1}/{3}}
\varepsilon^{{2}/{3}}\geqslant\tfrac{1}{2}\left(\tfrac{1}
{2}\Phi(x)\right)^{{1}/{3} }
\varepsilon^{{2}/{3}}>\tfrac{1}{3}\Phi(x)^{{1}/{3} }
\varepsilon^{{2}/{3}},
$$
so we obtain the statement of the theorem.
\qed

\PR{Proof of Theorem \ref{T:GlobalEquivalence}}
 Note first that Theorems \ref{T:ExplicitBounds},
\ref{T:LocalExample} have analogous
statements for the infimum, instead of the supremum, which clearly
can be derived from these theorems.

We have $\|\{f,g\}\|>0$, since otherwise every point in
$ M $ is a degenerate
critical point of $ \{ f,g \} $. Then for any $ 1 \leqslant k
\leqslant N $ we have $ \{ f,g \} ( x_k ) \neq 0 $, therefore in
particular $ x_k $ is not a critical point for each
of the functions $ f,g $. Therefore, we can apply
Theorem \ref{T:ExplicitBounds}, together with the remark
at the
beginning of the proof, to obtain the inequality
$$
\limsup_{ \varepsilon \rightarrow 0 }
\frac{\Upsilon_{f,g}(\varepsilon)}{\varepsilon^{{2}/{3}} }
\leqslant 6 \big| \Phi ( x_k )\big|^{{1}/{3}}.
$$
This is true for
any $ 1 \leqslant k \leqslant N $, so we obtain the desired upper bound.

Let us prove the lower bound.
For any $ 1 \leqslant k \leqslant N $, take a neighborhood
\hbox{$x_k\in U_{k}\subset M$}, such that
$|\{f,g\}(y)|<|\{
  f,g \} ( x_k ) |,$ for every $ y \in \overline{ U_k } \setminus \{ x_k \} $.
Then Theorem \ref{T:LocalExample} guarantees that there exist
neighborhoods $ x_k \in V_k \subset U_k $, such that for any $
\varepsilon $ small enough there exist functions $ F_k,G_k : M
\rightarrow \mathbb{R} $ satisfying
\begin{gather*}
\| F_k -f \| \leqslant
\varepsilon\,,\quad\|G_k-g\|\leqslant\varepsilon\,,\\
\{F_k,G_k\}(y)\leqslant\big\|\{f,g\}\big\|-\tfrac{1}{3}
\big|\Phi(x_k)\big|^{
{1}/{3} }\varepsilon^{{2}/{3}} ,\quad\forall y\in U_k\,,
\end{gather*}
and such that $ F_k = f, G = g_k $ on $ M \setminus V_k $.

\enlargethispage{2pt}
Define $F,G:M\rightarrow\mathbb{R}$ as $F=F_k$, $G=G_k$ on
each of $U_k$ and $F=f$, \hbox{$G=g$} on $ M \setminus
\bigcup_{k=1}^{N} U_k $.
Then on the union $ \bigcup_{k=1}^{N} U_k
$ we clearly will have $ | \{ F,G \} | \leqslant \| \{f,g\} \| -
\frac{1}{3}C\varepsilon^{{2}/{3}}$, and for the set $ M
\setminus \bigcup_{k=1}^{N} U_k $ we have $ \max_{ M \setminus
\bigcup_{k=1}^{N} U_k } | \{ F,G \} |  = \max_{ M \setminus
\bigcup_{k=1}^{N} U_k } | \{ f,g \} | < \| \{ f,g \} \| $, and
does not depend on $ \varepsilon $. Therefore, for small $
\varepsilon$ we have
$$
\big\|\{F,G\}\big\|\leqslant\big\|\{f,g\}\big\| -
\tfrac{1}{3}C\varepsilon^{{2}/{3}}.
$$
This example of $ F,G $ shows that
$$
\tfrac{1}{3}C\varepsilon^{{2}/{3} } \leqslant
\Upsilon_{f,g}(\varepsilon)\,.\eqno{\scriptstyle\square}
$$

\PR{Proof of Theorem \ref{T:UnicontinPrelim}}
 First of all, consider the case when
$$ \| G_1 \|_{U,1} = \| G_2
 \|_{U,1}= 0\,.
$$
In this case, we clearly have $ dG_1 = dG_2 = 0 $
 on $ U $, hence $$ \{ F_1, G_1 \} =  \{ F_2, G_2 \} = 0 $$ on $ U
 $, and then the desired inequality
 $$
\inf_{y,z\in U}\big|\{F_1,G_1\}(y)-\{F_2,G_2\} (z)\big|
\leqslant C\varepsilon\max\bigl(1,\|G_1\|_{U,1},\|G_2\|_
{U,1}\big)
$$
is satisfied for any choice of $ C > 0 $.

\goodbreak
We are left with the case of
$$
\max\bigl( \| G_1 \|_{U,1} , \| G_2
 \|_{U,1}\big) > 0\,.
$$
Denote by $ \Phi_{G_1}^{t},\Phi_{G_2}^{t} $ the Hamiltonian flows
corresponding to the Hamiltonians $ G_1,G_2 $. Take some open
subset $ V \subset U $, such that the closure $ \overline{V}
\subset U $.
Clearly there exists a constant $ c_1 = c_1 (U,V) $,
such that for any $$ 0 < t \leqslant T := \frac{c_1}{ \max ( \|
G_1 \|_{U,1} , \| G_2 \|_{U,1} ) }\,,
$$
and for any $ y \in V $, we
have that $ \Phi_{G_1}^{t}(y),\Phi_{G_2}^{t}(y) \in U $. Take some
$ \delta>0$, and assume that we have
$$
\inf_{y,z\in U}\bigl| \{
F_1 , G_1 \} (y) - \{ F_2 , G_2 \} (z)
\big| > \delta\,.
$$
Then one of the following holds:

either

(a) $ \inf_{z \in U } \{ F_2 , G_2 \} (z) - \sup_{y \in U} \{ F_1 ,
G_1 \} (y) > \delta $,

or

(b)  $ \inf_{y \in U} \{ F_1 , G_1 \} (y) - \sup_{z \in U} \{ F_2 ,
G_2 \} (z) > \delta $.

Assume for instance, that (a) holds.
Denote $K=\sup_{ y \in U } \{ F_1 ,G_1 \} $.
Fix any \hbox{$y\in V$}, and denote $ K(t) := F_1
(\Phi_{G_1}^{t}(y))$, for $t\in[0,T]$.
Then for every $ t
\in [0,T] $ we have
$ K'(t) = \{ F_1 ,G_1 \} (\Phi_{G_1}^{t}(y))
\leqslant K $, hence  for every $ t \in [0,T] $ we have
\hbox{$K(t)-K(0)$}
$=F_1(\Phi_{G_1}^{t}(y))-F_1(y)=\int_{0}^{t}K'(s)
ds \leqslant Kt $. Analogously, for any $ z \in V $, for any $ t
\in [0, T] $ we have
$ F_2(\Phi_{G_2}^{t}(z)) - F_2(z) \geqslant
(K + \delta)t $. Then we have
\begin{gather} \label{e:100}
f\big(\Phi_{G_1}^{t}(y)\big)-f(y)\leqslant Kt+2\varepsilon
\,,
\\
\label{e:101}
f\big(\Phi_{G_2}^{t}(z)\big) - f(z)
 \geqslant (K + \delta)t - 2 \varepsilon\,,
\end{gather}
for any $ y,z \in V $.
 Consider any point $ x \in V $ and for $ r,\alpha > 0 $ denote
$$
W=W_{r,\alpha}=B_{x}(r)\cap\big\{y\in M\mid f(x)<f(y)<f(x)
+\alpha\big\}\subset M\,,
$$
where $ B_{x}(r) $ is a ball of radius $ r $ centered at $ x
$, with respect to the metric $ \rho $. Then for small $ r,\alpha
$ we have $ W \subset V $. For any $ y,z \in W_{r,\alpha} $ we
have
\begin{equation} \label{e:102}
\big| f(y) - f(z)\big| < \alpha\,.
\end{equation}

From the inequalities (\ref{e:100}), (\ref{e:101}),
(\ref{e:102})
we conclude that for any $ y,z \in V $ we have
\begin{align*}
f\big(\Phi_{G_2}^{t}(z)\big)&\geqslant f(z)+(K+\delta)t- 2
\varepsilon > f(y) - \alpha + (K + \delta)t - 2 \varepsilon
\\
&\geqslant f\big(\Phi_{G_1}^{t}(y)\big) - Kt - 2
\varepsilon - \alpha + (K{+}\delta)t - 2 \varepsilon =f
\big(\Phi_{G_1}^{t}(y)\big)+\delta t-4\varepsilon-\alpha
\,.
\end{align*}
Therefore, if we assume that
\begin{equation} \label{e:103}
 \delta t \geqslant 4 \varepsilon +
\alpha\,,
\end{equation}
we get that $ f(\Phi_{G_2}^{t}(z)) > f(\Phi_{G_1}^{t}(y))
$ for any $y,z\in W$, therefore
$\Phi_{G_1}^{t}(W)\cap\Phi_{G_2}^{t}(W)$ $=\emptyset$,
hence the set $ W $ is displaced
by the map $\Phi_{G_1}^{-t}\circ\Phi_{G_2}^{t}$.
We have
the estimation
$$
\|\Phi_{G_1}^{-t}\circ\Phi_{G_2}^{t} \|_{Hof}
\leqslant2t \| G_2 - G_1 \|
 <2\varepsilon t $$
of the Hofer norm.
On the other hand, as a conclusion from
Lemma \ref{L:SymplecticEnergy} (see Definition
\ref{D:Norms2} for the notation used in the lemma), there
exists a constant $ c_2 =
c_2 (\rho,f,x)>0$, such that for small
$r,\alpha,{\alpha}/{r} $ we have $ e( W_{r, \alpha
} ) \geqslant c_2 r \alpha $. Therefore, we conclude that for $ t \in
[0,T]$, and small $r,\alpha,{\alpha}/{r}>0$, satisfying
(\ref{e:103}) we have
$$ c_2 r \alpha \leqslant e( W_{r, \alpha } )
\leqslant \| \Phi_{G_1}^{-t} \circ \Phi_{G_2}^{t} \|_{Hof} <
2\varepsilon t\,.
$$
Hence we conclude that given $ \delta, t, r, \alpha>0$,
satisfying (a), (\ref{e:103}), and $t\in[0,T]$, and
if $ r,
\alpha,{\alpha}/{r}$ are small enough, then we have $c_2 r
\alpha <2\varepsilon t $.
An analogous statement holds also for the
condition (b).
Therefore, we have

\begin{lem} \label{L:ConclusionUniCont}
 There exist constants $ c_1,c_2 > 0 $ such that for
any $\delta>0$, $0<t\leqslant\frac{c_1}{\max ( \| G_1
\|_{U,1} , \| G_2 \|_{U,1} ) } $, and small $ r, \alpha,
\frac{\alpha}{r} > 0 $, satisfying
\begin{gather*}
\inf_{y,z\in U}\bigl| \{ F_1 ,
G_1 \} (y) - \{ F_2 , G_2 \} (z)\big| > \delta\,, \\
\delta t\geqslant 4\varepsilon+\alpha\,,
\end{gather*}
we have $ c_2 r
\alpha <2\varepsilon t $.
\end{lem}

 Fix some small $ r = r_0 $, take $$ t =
\frac{c_1}{\max(\|G_1\|_{U,1},\|G_2\|_{U,1})+1 }\,,
$$
$\alpha=\frac{2\varepsilon t}{c_2r}$, and then take
$\delta=
\frac{4 \varepsilon + \alpha}{t} $. The value of $ r=r_0 $ is already
chosen to be small and fixed, and since
$t\leqslant c_1$, we have
$\alpha\leqslant\frac{2c_1}{c_2 r_0} \varepsilon $,
$\frac{\alpha}{r}=\frac{\alpha}{r_0}\leqslant\frac{2c_1}
{c_2 r_0^2}
\varepsilon $, that are small if $ \varepsilon $ is small.

Therefore, we can apply Lemma \ref{L:ConclusionUniCont},
and obtain
$$
\inf_{y,z\in U}\bigl|\{F_1,G_1\}(y)-\{F_2,G_2\}(z)\big|
\leqslant \delta\,.
$$
We have
$$
\delta=\frac{4 \varepsilon}{t} +
\frac{\alpha}{t}=\frac{4}{c_1}\varepsilon\max\big(\|G_1
\|_{U,1} ,\|G_2\|_{U,1}\big)+\left(\frac{4}{c_1}+\frac
{2}{c_2r_0}\right)\varepsilon\,.
$$
Therefore, denoting $C=\frac{8}{c_1}+\frac{4}{c_2r_0}$, we
obtain the statement of Theorem \ref{T:UnicontinPrelim}.
\qed

\PR{Proof of Theorem \ref{T:UniContin}}
 Consider any open $ U \subset M $, with compact closure $
\overline{U} \subset M $.
Take any $ n \in \mathbb{N} $ and apply
Theorem \ref{T:UnicontinPrelim} to the functions
$F_1=f_n$, $G_1 = g_n$, $F_2 = f$, $G_2 = g $.
We will get
\begin{multline*}
\inf_{y,z\in U}\big|\{f_n,g_n\}(y)-\{f,g\}(z)\big|\\
\leqslant C\cdot\max\bigl(  \| f_n -f \|_{U}  , \|g_n-g\|
_{U}\big)\cdot\max\big(1,\|g\|_{U,1},\|g_n\|_{U,1}\big)
\,.
\end{multline*}
Hence for some constant $ C' $ we have
\begin{multline*}
\inf_{y,z\in U}\big|\{f_n,g_n\}(y)-\{f,g\}(z)\big|\\
\leqslant C'\max\bigl(\|f_n-f\|_{U},\| g_n -g
\|_{U}\big)\|g_n\|_{U,1}+C'\max\big(\|f_n-f\|_{U},\|g_n-g
\|_{U}\big)\,.
\end{multline*}
Because of the assumptions of the theorem, the right-hand
side converges to $ 0 $, when $ n \rightarrow \infty $. On
the other hand, the sequence of functions $ \{ f_n , g_n \} $
uniformly converges to the function $ h $. Therefore, we conclude
that
$$
\inf_{y,z\in U}\big|h(y)-\{f,g\}(z)\big|=0\,.
$$
This holds for any open $ U \subset M $ with compact closure $
\overline{U} \subset M $. Then, because the functions $ h, \{ f,g
\} $ are continuous, we get that $ h(x) = \{ f,g \} (x) $ for any
point $ x \in M $.
\qed

\PR{Proof of Theorem \ref{T:FirstOrder}}
(a)
Since $ \lambda : C^{\infty}(X) \rightarrow C^{\infty}(X) $ is
a differential operator of the first order, there exists a vector
field $ v \in TX $ such that $ \lambda (f) = df(v) $.
There exists a positive $T=T(x,U)$, such that we have a
well-defined flow $
\Phi^{t}(x) $ of $ v $, for $ t \leqslant T $, and moreover $
\Phi^{t}(x) \in U $, for $ 0 \leqslant t \leqslant T $.
Assume that we are given $\varepsilon>0$ and a smooth
function $ F : M \rightarrow
 \mathbb{R} $, such that $ \| f-F \| \leqslant \varepsilon $.
 Denote $ K(t) = f( \Phi^{t}(x) )$, $ L(t) = F( \Phi^{t}(x) ) $.
Assume for a moment that we have some $\delta>0$ such that
$$
\lambda(F) \leqslant \lambda(f)(x) - \delta = K'(0) - \delta $$
on $ U $. Then  $ L'(t) \leqslant (K'(0)-\delta) $, hence $ L(t)
\leqslant L(0) + (K'(0)-\delta)t $, for $ t \leqslant T $. Because
of the assumption $ \| f-F \| \leqslant \varepsilon $, we have
\begin{align*}
K(t)&\leqslant L(t) + \varepsilon \leqslant
L(0)+\big(K'(0)-\delta\big)t+\varepsilon\leqslant\\
&\leqslant
K(0)+\big(K'(0)-\delta\big)t + 2\varepsilon\,,
\end{align*}
hence
$$ \delta t
\leqslant K(0) + K'(0)t - K(t) + 2\varepsilon\,.
$$
We have
$$
K(t)=K(0)+K'(0)t+\frac{1}{2}K''(0)t^{2}+\frac{1}{6}K'''
(0)t^{3}
+O(t^4)\,.
$$
On the other hand, $ K''(0) = 0 $, since the
function $ \lambda(f) $ attains its maximum at the point $ x $, and we
see that
\begin{align*}
&K'(t)=df\big(v(\Phi^{t}(x))\big)=\lambda(f)\big(\Phi^{t}
(x)
\big)\,,\\
&K''(t)=d (\lambda(f))\big(v( \Phi^{t}(x) )\big) =
\lambda^{2}(f)\big(\Phi^{t}(x)\big)\,,\\
&K'''(t)=d\big(\lambda^{2}(f)\big)\big(v(\Phi^{t}(x))\big
)=\lambda^{3}(f)\big(\Phi^{t}(x)\big)\,.
\end{align*}
Therefore,
$ \delta t \leqslant -\frac{1}{6}\lambda^{3}(f)(x)t^{3} + O( t^4 )
+ 2\varepsilon $, hence
$$
\delta\leqslant-\frac{1}{6}\lambda^{3}(f)(x)t^{2}+\frac
{2\varepsilon}{t} +
O(t^{3})\,,
$$
for every $ t \leqslant T $.
We substitute
$t=t_0=({-6\varepsilon}/{\lambda^{3}(f)(x)})^
{{1}/{3}}$
and we get
$$
\delta\leqslant 3\left(-\frac{1}{6}
\lambda^{3}(f)(x)\right)^{{1}/{3}}\varepsilon^{{2}/{3}}+O
(\varepsilon)=\left(\frac{9}{2}\right)^{{1}/{3}}\big(-
\lambda^{3}(f)
\big)^{{1}/{3}}\varepsilon^{{2}/{3}}+O(\varepsilon)\,.
$$
Note that $ t_0 < T $ , when $ \varepsilon $ is small.
 This observation leads to the desired result.

\smallskip
(b)
Since $B(\,\cdot\,{,}\,\cdot\,)$ is not anti-symmetric,
there
exists some $ h \in C^{\infty}(X) $ such that $ B(h,h) $ is a
non-zero function.
Take any smooth nondecreasing function $
\varphi:\mathbb{R}\rightarrow\mathbb{R}$, such that $
\varphi(t)=2n$, $\forall t\in[2n,2n+1]$, for every $ n \in
\mathbb{Z} $. Define $ f_{n}(x) = \frac{1}{n} \varphi ( n h(x)) ,
g_{n}(x) = \frac{1}{n} \varphi ( n h(x) +1 ) $. It is easy to see,
that $f_n,g_n\rightarrow h$ uniformly, but
\begin{align*}
B(f_n,g_n)(x)&=B\left(\frac{1}{n}\varphi(nh(x)),\frac{1}
{n}
\varphi\big(nh(x)+1\big)\right)\\
&=\varphi'(nh(x))\varphi'\big( n
h(x)+1\big)B(h,h)=0\,,
\end{align*}
since $ \varphi'(t)\varphi'(t+1)=0 \forall
t \in \mathbb{R} $.
\qed

\PR{Proof of Theorem \ref{T:Higher order}}
  Denote $ h = \{ f,g \} $.
 The proof goes similarly to that of Theorem
\ref{T:ExplicitBounds}.
We will use the notation in the Definition
\ref{D:Norms2},
introduced in the proof of Theorem \ref{T:ExplicitBounds}.

Instead of inequality (\ref{e:25}) we will have
\begin{align*}
\delta&\geqslant\frac{\|h\|_{U,2l}}{2l+1}\frac{(r+t\|X_{g}
  \|_{U})^{2l+1} }{ t \| X_{g} \|_{U} }
+\frac{2\varepsilon}{t}+\frac{\alpha}{t}\\
&=
\frac{\| h \|_{U,2l}}{2l+1} \frac{ (r + t \| X_{g}
  \|_{U})^{2l+1} }{ t \| X_{g} \|_{U} } + \frac{2\varepsilon}{t} +
\frac{2\varepsilon}{Cr}\,.
\end{align*}

Our choice of $ t,r $ will be of the form
$t={P\varepsilon^{\frac{1}{2l+1}}}/{\|X_{g}\|_{U}}$,
  $ r = P \varepsilon^{ \frac{1}{2l+1} } $, for some $ P > 0 $.
Then we have
\begin{multline*}
\frac{\| h \|_{U,2l}}{2l+1} \frac{ (r + t \| X_{g}
  \|_{U})^{2l+1} }{ t \| X_{g} \|_{U} } + \frac{2\varepsilon}{t} +
\frac{2\varepsilon}{Cr}\\
=\left(\frac{2^{2l+1}}{2l+1}\|h\|_{U,2l}P^{2l}+2\left(\|X
_{g}\|_{U}+\frac{1}{C}\right)\frac{1}{P}\right)
\varepsilon^{
\frac{2l}{2l+1} }\,.
\end{multline*}
We fix $P$, that minimizes the expression
$$
\frac{2^{2l+1}}{2l+1} \| h \|_{U,2l}
P^{2l}+2\left(\|X_{g}\|_{U}+\frac{1}{C}\right)\frac{1}{P}
\,.
$$
The corresponding value of $P$ does not depend on
$\varepsilon$.
Then we take $ \varepsilon $ small enough, such that the
assumptions of Lemma \ref{L:SymplecticEnergy} are
satisfied,
  and we obtain
$$
\frac{4l+2}l\left(\frac l{2l+1}\right)^{\frac 1{2l+1}}
\bigg(\|h\|_{U,2l}\left(\|X_{g}\|_{U}+\frac{1}
{C} \right)^{2l} \bigg)^{ \frac{1}{2l+1} }
\varepsilon^{ \frac{2l}{2l+1}}.
$$
Then, by the same arguments
  as in Theorem \ref{T:ExplicitBounds} we arrive at
$$
\limsup_{ \varepsilon \rightarrow 0 }
\frac{\Upsilon_{f,g}^{+}(\varepsilon)}{
\varepsilon^{\frac{2l}{2l+1}}}\leqslant\frac{8l+4}l\left(
\frac l{4l+2}\right)^{\frac1{2l+1}}\left(\frac{1}{2l!}
\max_{\theta} P_{2l}( \theta ) \right)^{ \frac{1}{2l+1} }
  , $$
where $ P_{2l}( \theta ) $ equals
$$
-\big\{\dots\{\{h,\cos(\theta)f+\sin(\theta)g\},\cos(
\theta
)f + \sin( \theta )g\}
,\ldots,\cos(\theta)f+\sin(\theta)g\big\}(x)\,,
  $$
when the Poisson bracket is taken $ 2l $ times.
Note that $ P_{2l} $ is a non-negative
 trigonometric polynomial of degree $ \leqslant 2l $.

 \begin{lem} \label{L:TrigPolNorm^2}
 There exists a complex trigonometric polynomial $ Q(\theta) $ of degree $ \leqslant l $, such that
 $$ P_{2l}(\theta) = |Q(\theta)|^2 .$$
 \end{lem}

 \PR{Proof of Lemma \ref{L:TrigPolNorm^2}}
Let us remark, that along the proof we will only use the
fact that $ P_{2l}(\theta) $ is non-negative.

  Denoting $ z = \cos( \theta ) + i \sin( \theta ) $, the
  trigonometric polynomial $ P_{2l}(\theta) $ can be written as a
  polynomial of $ z,{1}/{z} $, and there exists a complex
  polynomial $ T \in \mathbb{C}[z] $,
  such that $$ P_{2l}(\theta) = \frac{1}{z^r} T(z) ,$$ and $T(0)
  \neq 0 $.
Since $ P_{2l}(\theta) $ is a real number for any $ \theta \in
  \mathbb{R} $ , then for any $ z \in \mathbb{C} $, $ |z|
  = 1 $, we have that $  \frac{1}{z^r} T(z) \in \mathbb{R} $, hence
  $$ \frac{1}{z^r} T(z) =  \overline{ \frac{1}{z^r} T(z) } .$$
Assume that $ T(z) = c \prod_{k=0}^{m} (z - \alpha_k) $.
Since $ T(0) \neq 0 $, we have that
\hbox{$\alpha_1,\alpha_2,\ldots,\alpha_m\neq 0$.}
\goodbreak

   Then for any $ z \in \mathbb{C} $ with $ |z| = 1 $, we have
\begin{align*}
     \overline{ \frac{1}{z^r} T(z) }
&=\frac{1}{\overline{z}^r}\overline{c}\prod_{k=0}^{m}(
\overline{z}-\overline{\alpha_k})=z^r\overline{c}\prod_{k
=0}^{m}\left(\frac{1}{z}-\overline{\alpha_k}\right)
\\
&=\frac{(-1)^{m}\overline{c}}{\prod_{k=0}^{m}\overline{
\alpha_k}}z^{r-m}\prod_{k=0}^{m}\left(z-\frac{1}{
     \overline{\alpha_k}}\right).
\end{align*}
Denote
$c'=\frac{(-1)^{m}\overline{c}}{\prod_{k=0}^{m}\overline{
\alpha_k}
}$.
Then
$$
\frac{1}{z^r}T(z)=\overline{\frac{1}{z^r}T(z)}=c'z^{r-m}
\prod_{k=0}^{m}\left( z - \frac{1}{
  \overline{\alpha_k}}\right),
$$
and hence
$$
z^mT(z)-c'z^{2r}\prod_{k=0}^{m}\left(z-\frac{1}{
  \overline{\alpha_k}}\right)=0\,,
$$
for any $ z \in \mathbb{C} $ with $ |z| = 1
  $.
Since a non-zero polynomial must have a finite number of
  roots, we must have an identity
$$
z^mT(z)=cz^m\prod_{k=0}^{m}(z-\alpha_k)\equiv c'z^{2r}
\prod_{k=0}^{m}\left( z - \frac{1}{
  \overline{\alpha_k}}\right) ,
$$
  as polynomials. Hence the list
  $$ \frac{1}{\overline{\alpha_1}} , \frac{1}{\overline{\alpha_2}} ,\ldots,
  \frac{1}{\overline{\alpha_m}} $$
is a permutation of
$$ \alpha_1 , \alpha_2 ,\ldots,
  \alpha_m\,.
$$
Moreover, if some $\alpha_j$ satisfies $|\alpha_j | = 1 $,
then its multiplicity as a root of the poly\-nomial
$T(z)$, is
  even. Indeed, write $ \alpha_j = e^{ i \theta_j } $, $\theta_j \in \mathbb{R} $, and consider the limit
\begin{align*}
  \lim_{ \tau \rightarrow 0}\frac{P_{2l}(\theta_j +
  \tau) }{ P_{2l}(\theta_j - \tau) }
& =
  \lim_{ \tau \rightarrow 0 }
  \frac{ \alpha_j^r e^{ -i r  \tau } }{ \alpha_j^r e^{ i r  \tau } }
  \lim_{ \tau \rightarrow 0 } \frac{ T( \alpha_j e^{ i  \tau } )}{ T( \alpha_j e^{- i \tau }
  ) }
  =
  \lim_{ \tau \rightarrow 0 } \frac{ T( \alpha_j e^{ i  \tau } )}{ T( \alpha_j e^{- i \tau
  }) }
\\
&=
  \lim_{ \tau \rightarrow 0 }
  \frac{ c \prod_{k=0}^{m} ( \alpha_j e^{ i  \tau } - \alpha_k )
  }{ c \prod_{k=0}^{m} ( \alpha_j e^{ - i  \tau } - \alpha_k ) }
  =
  \prod_{k=0}^{m}
\lim_{\tau\rightarrow 0}\frac{e^{i\tau}-\alpha_k\alpha_j^
{-1}}{e^{-i\tau}-\alpha_k\alpha_j^{-1}}\,.
\end{align*}
  We have that each of the terms $ \lim_{ \tau \rightarrow 0 }
  \frac{ e^{ i  \tau } - \alpha_k \alpha_j^{-1} }{ e^{ - i  \tau } - \alpha_k \alpha_j^{-1}
}$ equals $1$ if $\alpha_k\neq\alpha_j$, and $-1$ if
  $ \alpha_k = \alpha_j $.
Therefore, the limit equals $1$ if the multiplicity of
$ \alpha_k $ is even,
and $-1$ if the multiplicity of $\alpha_k$ is odd.
On the
  other hand, the limit $
  \lim_{ \tau \rightarrow 0 } \frac{ P_{2l}(\theta_j +
  \tau) }{ P_{2l}(\theta_j - \tau) } $ must be non-negative, because
  the trigonometric polynomial $ P_{2l} $ is non-negative. This
  proves, that the multiplicity $ \alpha_j $ is even.

As a conclusion, we obtain that the list of roots $$ \alpha_1 , \alpha_2 ,\ldots,
  \alpha_m $$
splits into pairs $\beta_j,\gamma_j$, $j=1,2,\ldots,s$,
such that $\gamma_j={1}/{\overline{\beta_j}}$, for every
$1\leqslant j \leqslant s
$, where $2s=m$.
Denote
$$
q(z) := \prod_{k=0}^{s} (z - \beta_k)\,.
$$
Then for $z = \cos( \theta ) + i \sin( \theta ) $, we have
\begin{align*}
q(z) \overline{q(z)}
&=\prod_{k=0}^{s}(z-\beta_k)\overline{\prod_{k=0}^{s} (z -
    \beta_k) }
\\
&=\prod_{k=0}^{s}(z-\beta_k)\prod_{k=0}^{s}\left(\frac{1}
{z}-\overline{\beta_k}\right)=\prod_{k=0}^{s}(z-\beta_k)
\prod_{k=0}^{s}\left(\frac{1}{z} - \overline{\beta_k}\right)
 \\
&
=\bigg((-1)^s\prod_{k=0}^{s}\overline{\beta_k}\bigg)\frac
{1}{z^s} \prod_{k=0}^{s}
(z - \beta_k) \prod_{k=0}^{s} (z - \gamma_k)
\\
&=\bigg((-1)^s\prod_{k=0}^{s}\overline{\beta_k}\bigg)
\frac{1}{z^s} \prod_{k=0}^{m}
  (z - \alpha_k)
  = \frac{ (-1)^s \prod_{k=0}^{s} \overline{\beta_k} }{c} \frac{1}{z^s} T(z)
  \\
&
  = \frac{ (-1)^s \prod_{k=0}^{s} \overline{\beta_k} }{c} z^{r-s}
  P_{2l}(\theta)
\,.
\end{align*}
  Denote $ c'' := \frac{ (-1)^s \prod_{k=0}^{s} \overline{\beta_k}}{c}
  $.
Then since we have that $ q(z) \overline{q(z)} ,
  P_{2l}(\theta) > 0 $ for any $ \theta $, except, may be, a
  finite number of values, therefore $ c'' z^{r-s} $ is a positive
  real number, for any $ z \in \mathbb{C} $, $ |z| = 1 $, possibly except a
  finite number of values. As a consequence, we have that $ r = s
  $, and $ c'' $ is a positive real number. Hence
  $$ q(z) \overline{q(z)}
  = c''  P_{2l}(\theta)
\,,
$$
and if we denote $Q(\theta) := \frac{1}{ \sqrt{ c'' } } q(
  \cos(\theta) + i \sin(\theta) ) $, we obtain
$$
|Q(\theta)|^2 = Q(\theta)\overline{Q(\theta)}= P_{2l}(
\theta)\,.\eqno{\scriptstyle\square}
$$

 \begin{lem} \label{L:TrigoIneq}
 $$ \max_{\theta} P_{2l}(\theta) \leqslant
\frac{2l+1}{2\pi}\int_{0}^{2\pi}P_{2l}(\theta)d\theta\,.
$$
 \end{lem}

 \PR{Proof of Lemma \ref{L:TrigoIneq}}
Because of Lemma \ref{L:TrigPolNorm^2} there exists a
complex
  trigonometric polynomial $ Q(\theta) $ of degree $ \leqslant l $, such that
  $ P_{2l}(\theta) =
  |Q(\theta)|^{2} $. Denote by $ a_{-l},a_{-l+1},\ldots,a_{l} $ the
  Fourier coefficients of $ Q(\theta) $.
Then by Holder
  inequality, for any $ \phi $ we have
\begin{align*}
P_{2l}(\phi)&=|Q(\phi)|^2=\big|a_{-l}e^{-il\phi}+a_{-l+1}
e^{- i (l-1) \phi} + \dots
+a_{l}e^{il\phi}\big|^2\\
&\leqslant\big(|a_{-l}|^2 +|a_{-l+1}|^2+\dots+a_{l}|^2
\big)(2l+1)=\frac{2l+1}{2\pi}\int_{0}^{2\pi}|Q(\theta)|^2
d\theta\\
&=\frac{2l+1}{2\pi}\int_{0}^{2\pi}|Q(\theta)|^2 d\theta =
\frac{2l+1}{2\pi}\int_{0}^{2\pi}P_{2l}(\theta)d\theta\,.
\tag*{$\scriptstyle\square$}
\end{align*}

Assume that
 $P_{2l}(\theta)=\sum_{k=0}^{2l}c_{k}\cos(\theta)^{2l-k} \sin(\theta)^{k}
$.
We have
$$
\int_{0}^{2\pi} \cos(\theta)^{2l-k} \sin(\theta)^{k}
d\theta=0\,,
$$
when $k$ is odd, and
$$
\int_{0}^{2\pi} \cos(\theta)^{2l-k}
 \sin(\theta)^{k} d\theta
=2\mathbf{B}\left(\frac{k}{2}+\frac{1}{2},l-\frac{k}{2}+
\frac{1}{2}\right),
$$
for even $ k $, where $ \mathbf{B}(x,y) $ is
the beta-function. It is easy to see that for any\break
$0\leqslant k\leqslant 2l$, we have that $c_k$ equals
the sum of terms of the form\break $ - \{\dots  \{ \{ h , f_1 \} , f_2 \}
,\dots  \} , f_{2l} \} (x) $, when each of the functions $ f_j $ is
one of $ f,g $, while the function $ f $ occurs $ 2l-k $ times,
and $ g $ occurs $ k $ times.
Since $ h $ has
multiplicity at least $ 2l $ at the point $ x $, all
these terms are equal. Indeed, for any $ 1 \leqslant m < 2l $,
denoting
$$
H=\big\{\{\dots  \{ \{ h , f_1 \} , f_2 \} ,\dots\} ,
f_{m-1}\big\}\,,
$$
we have
$$
\big\{ \{ H, f_{m} \} , f_{m+1}\big \} =
\big\{\{H,f_{m+1}\},f_{m}\big\}+\big\{H,\{f_{m},f_{m+1} \}
\big\}\,,
$$
hence
\begin{align*}
&\big\{\dots\{\{\{\dots\{\{h,f_1\},f_2\},\dots\},
f_{m-1} \} , f_{m} \} ,\ldots, f_{2l}\big\}\\
&\qquad=\big\{\dots  \{ \{ \{
H,f_{m}\},f_{m+1}\},f_{m+2}\},\ldots,f_{2l}\big\}\\
&\qquad=\big\{\dots\{\{H,f_{m+1}\},f_{m}\},\ldots,f_{2l}
\big\}+\big\{\dots
\{\{H,\{f_{m},f_{m+1}\}\},f_{m+2}\}\dots,f_{2l}\big\}
\\
&\qquad=\big\{\dots\{\{\{\{\{\dots\{\{h,f_1\},f_2\},\dots
\},f_{m-1}
\},f_{m+1}\},f_{m}\} , f_{m+2} \} ,\ldots, f_{2l}\big\}\\
&\qquad\qquad+\big\{\dots\{\{\{\{\dots\{\{h,f_1\},f_2
\},\dots\},f_{m-1}\},\{f_{m} ,
f_{m+1} \} \} , f_{m+2} \} ,\ldots, f_{2l}\big\}\,,
\end{align*}
 and
$$
\big\{\dots\{\{\{\{\dots\{\{h,f_1\},f_2\},\dots\big\},f_
{m-1}\},\{
f_{m},f_{m+1}\}\},f_{m+2}\} ,\ldots, f_{2l}\big\}(x)= 0\,,
$$
since we have applied the Poisson bracket $ 2l-1 $ times, starting
with the function $ h $, and $ h $ has multiplicity $ 2l $ at $ x
$.
Therefore, we have that
$$
c_k = \left( {\begin{array}{*{20}c}
2l \\ k \\ \end{array}} \right) H_{k}(x)  = \frac{1}{
\mathbf{B}(k,2l-k)} H_{k}(x)\,,$$
where
$$
H_{k}=-\big\{\dots\{\{h,f\},f\},\dots\big\},f\},g\},g
\} ,\ldots, g \}\,,
$$
when $ f $ appears $ 2l-k $ times, and $ g $
appears $ k $ times.
From all these observations we have
$$
\int_{0}^{2\pi} P_{2l}(\theta)
d\theta=2\sum_{m=0}^{l}\frac{\mathbf{B}\big(m+\frac{1}{2},
l-m+\frac{1}{2}\big)}{\mathbf{B}(2m,2l-2m)}H_{2m}(x)\,.
$$
Using
the identities, concerning the $ \mathbf{B} $ and $
\mathbf{\Gamma} $-functions, one can check that
$$
\frac{\mathbf{B}\big(m+\frac{1}{2},l-m+\frac{1}{2}\big)}{
\mathbf{B}(2m,2l-2m)} = \left( {\begin{array}{*{20}c} l \\ m \\
\end{array}} \right) .$$ Again, because $ h $ has multiplicity $ 2l $ at $ x $
, we have
$$
\sum_{m=0}^{l} \left( {\begin{array}{*{20}c} l \\ m \\
\end{array}}\right)H_{2m}(x)=-\mathcal{D}^{l}\big(\{f,g\}
\big)(x)\,.
$$
Summarizing the above considerations, we get that
\begin{align*}
\limsup_{\varepsilon\rightarrow 0}
\frac{\Upsilon_{f,g}^{+}(\varepsilon)}{\varepsilon^{\frac
{2l}{2l+1}}}&\leqslant-\frac{8l+4}l\left(\frac{l}{4l+2}
\right)^{\frac 1{2l+1}}\left(\frac{2l+1}{\pi}\right)^{
\frac1{2l+1}}\left(\frac1{2l!}\mathcal{D}^{l}\big(\{f,g\}
\big)(x)\right)^{
\frac{1}{2l+1} }\\
&\leqslant-9\left(\frac{1}{2l!}\mathcal{D}^{l}\big(\{f,g\}
\big)(x)\right)^{ \frac{1}{2l+1} }
\, .\tag*{$\scriptstyle\square$}
\end{align*}

\section{Non-locality} \label{S:Int and Loc}
 On first sight it seems that the statement of Theorem \ref{T:ExplicitBounds} is local, in the sense that if the
 Poisson bracket $ \{ f,g \} $ attains its maximum at the point $ x
 \in M $, then for any two sequences $$ f_1 , f_2 ,\ldots, g_1 , g_2
 ,\ldots  \in C^{\infty}(M)\,,
$$
such that $ \| f_n -f \| \rightarrow 0
 $, $ \| g_n -g \| \rightarrow 0 $, there exists a sequence $ x_n
 \rightarrow x $,  such that
 $$ \liminf_{ n \rightarrow \infty } \{ f_n , g_n \} (x_n) \geqslant
 \max\{f,g\}\,.
$$
In fact, we cannot conclude that, since the flow
 $ \Phi_{G}^{t} $ can be very fast, and during a small time can
 exit a neighborhood of $ x $. Actually, the locality does not hold for any dimension $ n>2 $.
 For dimension $ 2 $ the locality was proved by Zapolsky  \cite{Z}.

 On the other hand, the condition of existence of the flow $
 \Phi_{G}^{t} $ for all time $ t $ is essential, as we will see in
 the example below.

 The examples that reflect both of the remarks above are based on
 the example of Polterovich, mentioned in
 Example \ref{E:PolterovichExample}.

\begin{exmp} \label{E:Gcondition}
Consider the manifold
$$
M=\big\{(x,y,z,u)\in\mathbb{R}^{4}\mid 1<z<1\big\}\subset
 \mathbb{R}^{4} ,
$$
endowed with the standard symplectic form
$\omega=dx\wedge dy +
 dz \wedge du $.
Let $ \chi (t) := \sqrt{2t+2 }$, $t\in(-1,+\infty)$.
Then $ \chi(t) \chi'(t)=
 1 $.
Consider the functions
$$
f(x,y,z,y)=x\,,\quad g(x,y,z,u)=y\,,
$$
and define
 $$
f_{n}(x,y,z,u)=x+\frac{\chi(z)}{\sqrt{n}}\cos(nu)\,,\quad
g_n(x,y,z,u)=y-\frac{\chi(z)}{\sqrt{n}}\sin(nu)\,,
$$
for $ n=1,2,3,\dots$.
Then $ f_n \rightarrow f $,
 $ g_n \rightarrow g $ uniformly on $ M $. However, we have
 $ \{ f,g \} \equiv 1 $, but $ \{ f_n,g_n \} \equiv 0 $ for
 every $ n $, so rigidity does not hold in its weakest sense.

 The reason is that the flows $ \Phi_{g_n}^{t} $ are not defined
 for arbitrary time $ t $.
 \end{exmp}

 As a corollary of Example \ref{E:Gcondition}, we derive the non-locality of Theorem \ref{T:ExplicitBounds}.
 We already see the non-locality in Example \ref{E:Gcondition}, however, $ g_n $ does not belong to
 $ \mathcal{H}^{b}(M,\omega) $. One can fix this problem by the following truncation of
the functions.

 \begin{exmp} \label{E:NonLocal}
Consider the manifold $M$ and functions
$$
f,g,f_n,g_n : M \rightarrow \mathbb{R}\,,
$$
$n=1,2,\ldots$, as in the previous Example
  \ref{E:Gcondition}. Take a smooth function $ \psi : \mathbb{R} \rightarrow \mathbb{R}
 $, such that $\psi(x)=1$ for $|x|\leqslant{1}/{4}$,
$\psi(x)=0$ for $|x|\geqslant{1}/{3}$, and
$x\psi'(x)\leqslant 0 $ for all $ x $.
Then define $ \varphi : \mathbb{R}^{4} \rightarrow \mathbb{R} $ by
 $ \varphi(x,y,z,u) = \psi(x)\psi(y)\psi(z)\psi(u) $. Then $ x\varphi_{x},y\varphi_{y}
 \leqslant 0$. Denote
\begin{alignat*}{2}
&F(p)=f(p)\varphi(p)\,, &\quad&G(p)=g(p)\varphi(p)\,,\\
&F_n(p)=f_n(p)\varphi(p)\,,&\quad&G_n(p)=g_n(p)\varphi(p)
\,,
\end{alignat*}
for $ n=1,2,3,\ldots$, and $ p \in M $.
 Then $ F,G,F_n,G_n $ are all compactly supported.
 We have
\begin{align*}
\{F,G\}&=  \{ f\varphi,g\varphi \} =  \varphi^2
 + \varphi y \{ x ,\varphi \} + \varphi x \{ \varphi , y\}
&=\varphi^2+\varphi y\varphi_{y}+\varphi x\varphi_{x}
 \leqslant \varphi^2 \leqslant 1
\end{align*}
at every point, and $ \{ F,G \}
 = 1 $ in the cube $K:=\{|x|,|y|,|z|,|u|<{1}/{4}\}$.
 However, for every $ p \in K $, we have the equality $ F_n=f_n, G_n =
 g_n $, hence $ \{ F_n , G_n \} = 0 $ in $ K $. This reflects the
 non-locality. Note that
$$
\operatorname{supp}F,G,F_n,G_n\subset\left\{|x|,|y|,|z|,
|u|\leqslant
\tfrac{1}{3}\right\}.
$$
 Hence non-locality holds for any symplectic manifold of dimension
 $ 4 $, because of the existence of a Darboux chart on $ M $, and re-scaling of $ F,G,F_n,G_n $, in order that their supports
 be contained in this chart. Surely
 this is true in any dimension of $ M $, since one can provide a similar example for any even dimension
 bigger than $ 4 $.
 \end{exmp}

\subsection*{Acknowledgments.}
I would like to thank my supervisor Paul Biran for the help and
attention he gave me. I thank Leonid Polterovich, Michael Entov,
Michail Sodin, Frol Zapolsky, Oleg Khasanov and Egor Sheluhin for
helpful discussions. Also I would like to thank Alexander
Bykhovsky and Alexander Sodin for improving the style of the
paper.
And I thank Sobhan Seyfaddini for important remarks.

\bigskip\bigskip

\noindent{\sc Lev Buhovski},
The Mathematical Sciences Research Institute,
Berkeley, CA 94720-5070, USA
\hfill{\tt levbuh@gmail.com}

\bigskip
\rightline{Received: October 4, 2008}
\rightline{Revision: January 29, 2009}
\rightline{Accepted: February 4, 2009}


\small
\begin{thebibliography}{XXX}
\itemsep=0pt
\parskip=1pt

\bibitem[CV]{CV}{\sc F.\ Cardin, C.\ Viterbo}, {Commuting
Hamiltonians and Hamilton--Jacobi multi-time equations},
Duke Math.\ J. 144 (2008), 235--284.

\bibitem[EP1]{EP-1}{\sc M.\ Entov, L.\ Polterovich},
{$C^0$-rigidity of Poisson brackets},
 preprint; http://arxiv.org/abs/0712.2913

\bibitem[EP2]{EP-2}{\sc M.\ Entov, L.\ Polterovich},
{$C^0$-rigidity of the double Poisson
 bracket}, preprint; http://arxiv.org/abs/0807.4275

\bibitem[EPZ]{EPZ}{\sc M.\ Entov, L.\ Polterovich, F.\
Zapolsky}, {Quasi-morphisms
    and the Poisson bracket}, Pure and Applied Mathematics
    Quarterly {3:4} (2007), 1037--1055.

\bibitem[H]{H}{\sc V.\ Humili\`{e}re}, {Hamiltonian
pseudo-representations}, preprint; math/0703335, 2007.

\bibitem[MS]{MS}{\sc D.\ McDuff, D.\ Salamon},
{Introduction to Symplectic Topology}, 2-nd edition,
    Oxford University Press, New York, 1998.

\bibitem[Z]{Z}{\sc F.\ Zapolsky}, {Quasi-states and the
Poisson bracket on surfaces},
J.\ of Modern Dynamics {1:3} (2007), 465--475.

\end{thebibliography}
\end{document}